\documentclass[aop,preprint]{imsart}

\RequirePackage[numbers]{natbib}
\usepackage{amsmath,amssymb,amsthm}

\usepackage[T1]{fontenc}
\usepackage{tikz}
\usetikzlibrary{arrows.meta,decorations.pathreplacing,calc}

\newcommand{\Exf}[1]{\mathbb{E}\left[#1\right]}
\newcommand{\dd}{\mathrm{d}}
\newcommand{\Hb}{\mathbf{H_2}}
\newcommand{\w}{\mathbf{w}}
\newcommand{\e}{\mathrm{e}}
\newcommand{\gs}{\mathcal{G}}
\newcommand{\ds}{\mathcal{D}}

\startlocaldefs
\numberwithin{equation}{section}

\theoremstyle{plain}
\newtheorem{theorem}{Theorem}
\newtheorem{lemma}{Lemma}[section]
\newtheorem{corollary}[lemma]{Corollary}
\newtheorem{proposition}[lemma]{Proposition}

\theoremstyle{definition}
\newtheorem{definition}[lemma]{Definition}

\newtheorem{remarkx}[lemma]{Remark}
\newenvironment{remark}[1][]{\begin{remarkx}[#1]\mbox{}
}{\hfill$\diamond$\end{remarkx}}
\endlocaldefs

\usepackage[hidelinks]{hyperref}
\begin{document}

\begin{frontmatter}
	\title{The optimal sub-Gaussian normalisation for randomised monotone functions}
	\runtitle{Randomised monotone functions}
	\begin{aug}
		\author[A]{\fnms{Thomas} \snm{Anton}\ead[label=e1]{tca2130@columbia.edu}}
		\author[B]{\fnms{Rabee} \snm{Tourky}\ead[label=e2]{rabee.tourky@anu.edu.au}}
		\address[A]{Department of Economics, Columbia University, New York, USA.
			\printead{e1}}
		\address[B]{Research School of Economics, Australian National University,
			Canberra, Australia.
			\printead{e2}}
	\end{aug}

	\begin{abstract}
		Let $\mathcal{M}$ denote the class of randomised monotone functions on
		$\mathbb{R}$ with values in $[0,1]$, and let
		$U_{\mathcal{M}}\colon \mathbb{R}_+\to \mathbb{R}_+$ be the minimal
		function for which
		\[
			\mathbb{P}\left\{ \sqrt{\eta_f}\, \sup_{t\in\mathbb{R}}
			\left| f_Z(t) - \Exf{f_Z(t)} \right|
			\ge \varepsilon\sqrt{U_{\mathcal{M}}(\eta_f)} \right\}
			\le 2\e^{-2\varepsilon^2}
		\]
		holds for every member $f_Z$ of $\mathcal{M}$ with finite effective sample size
		$\eta_f$ and every positive $\varepsilon$. We prove that for every
		$x> 1$,
		\[
			\left| \sqrt{U_{\mathcal{M}}(x)} - \sqrt{\log_4 x} \right|
			\le 2 \min\!\left\{ 1,\, \frac{2 \ln(\e + \ln x)}{\sqrt{\ln x}} \right\}\,.
		\]
		The optimal adjustment $\sqrt{U_{\mathcal{M}}(x)}$ matches
		$\frac{1}{\sqrt{2\ln 2}}\sqrt{\ln x}$ for all  $x>1$,
		with residuals bounded as above.
	\end{abstract}

	\begin{keyword}[class=MSC2020]
		\kwd{Primary 60E15}
		\kwd{ Secondary 05D40}
		\kwd{60F10}
		\kwd{60F17}
		\kwd{60G15}
	\end{keyword}

	\begin{keyword}
		\kwd{concentration inequalities}
		\kwd{randomised monotone functions}
		\kwd{empirical distribution functions}
		\kwd{sharp constants}
		\kwd{effective sample size}
	\end{keyword}

\end{frontmatter}

\section{Introduction}\label{sec:intro}

Massart~\cite{massart1990tight} establishes that for every empirical
distribution function $F_n \colon \mathbb{R} \to [0,1]$ based on $n$
independent and identically distributed observations,
\[
	\mathbb{P}\left\{ \sqrt{n} \sup_{t\in\mathbb{R}} \left|F_n(t) -
	\Exf{F_n(t)}\right| \ge \varepsilon \right\} \le
	2\e^{-2\varepsilon^2}, \qquad \text{for all } \varepsilon > 0 \,.
\]
The normalisation $\sqrt{n}$ is the same as the optimal distribution-free
normalisation implied by Hoeffding's inequality~\cite{hoeffding1963probability}
for a single threshold:
\[
	\mathbb{P}\left\{ \sqrt{n} \left|F_n(t) - \Exf{F_n(t)}\right|
	\ge \varepsilon \right\} \le 2\e^{-2\varepsilon^2}, \qquad
	\text{for all } \varepsilon > 0 \,.
\]

The Massart bound hinges on the precise structure of the randomised function:
\[
	F_n(t) = \frac{1}{n} \sum_{i=1}^n \mathbf{1}_{\{Z_i \le t\}}, \qquad
	\mathbb{P}\{Z_i \le t\}= \mathbb{P}\{Z_j \le t\},\qquad
	\text{for each } t\in \mathbb{R}\,.
\]
In this paper we determine the optimal sub-Gaussian normalisation for
the class of all randomised \mbox{$[0,1]$-valued} monotone functions.

Let
\[
	f \colon \mathcal{Z}_1 \times \cdots \times \mathcal{Z}_n \times \Theta
	\to \mathbb{R}
\]
be a bounded function. Let  $\mathcal{Z} = \prod_{i=1}^n \mathcal{Z}_i$
denote the product  of the measurable spaces $\mathcal{Z}_i$;
the set $\Theta$ is the index set.  For each $z \in \mathcal{Z}$ and $\theta \in \Theta$,
define the sections $f_\theta(z) := f(z, \theta)$ and
$f_z(\theta) := f(z, \theta)$. We suppose throughout that, for every $\theta \in \Theta$,
the section \mbox{$f_\theta \colon \mathcal{Z} \to \mathbb{R}$} is measurable.

For a random vector $Z = (Z_1, \dots, Z_n)$ with independent components
$Z_i \in \mathcal{Z}_i$, the mapping $f_Z \colon \Theta \to \mathbb{R}$ is a
\emph{randomised function}, and the pair $(f, Z)$ is a \emph{randomisation}.

Define $\eta_f \in (0, \infty]$ by
\[
	\eta_f := \left( \sum_{i=1}^n c_i(f)^2 \right)^{-1}, \qquad \text{where} \qquad
	c_i(f) := \sup_{(\theta, z_i^\ast, z) \in \Theta \times \mathcal{Z}_i
		\times \mathcal{Z}}
	\left| f_\theta(z_i^\ast, z_{-i}) - f_\theta(z) \right| \,.
\]
If $\eta_f$ is finite, then the McDiarmid--Hoeffding
inequality~\cite{mcdiarmid_1989} gives, at a single parameter $\theta$,
\begin{equation}\label{eq:macd}
	\mathbb{P}\left\{ \sqrt{\eta_f}
	\left| f_Z(\theta) - \Exf{f_Z(\theta)} \right| \ge \varepsilon \right\}
	\le 2\e^{-2\varepsilon^2}, \qquad
	\text{for all } \varepsilon > 0 \,.
\end{equation}
For the empirical distribution function $F_n$ we have $\eta_{F_n} = n$.
When $\eta_f=\infty$, the function $f_Z$ is deterministic for every $Z$.
Because $f$ is bounded, $\eta_f>0$. We call $\eta_f$ the
\emph{effective sample size} and take $\sqrt{\eta_f}$ as the
baseline distribution-free normalisation (Remark~\ref{rem:lattice}).

Let $\mathsf{M}$ be the set of all deterministic non-decreasing functions on
$\mathbb{R}$ with values in $[0,1]$. Write $\mathcal{M}$ for the class of all
randomisations $(f, Z)$, ranging over all finite products of measurable spaces,
that satisfy
\[
	f_z \in \mathsf{M}, \qquad  \text{for each } z \in \mathcal{Z}\,.
\]
For convenience, write $f_Z \in \mathcal{M}$ and call
$f_Z \colon \mathbb{R} \to [0,1]$ a \emph{randomised monotone function}.

Let $U_{\mathcal{M}} \colon \mathbb{R}_+ \to \mathbb{R}_+$ be the infimum of all
functions $U \colon \mathbb{R}_+ \to \mathbb{R}_+$ that satisfy
\begin{equation}\label{eq:admis}
	\mathbb{P}\left\{ \sqrt{\eta_f} \sup_{t \in \mathbb{R}}
	\left| f_Z(t) - \Exf{f_Z(t)} \right|
	\ge \varepsilon \sqrt{U(\eta_f)} \right\}
	\le 2\e^{-2\varepsilon^2}, \qquad \text{for all } \varepsilon > 0 \,,
\end{equation}
for every $f_Z\in \mathcal{M}$ satisfying $\eta_f<\infty$.
The function $U_{\mathcal{M}}$ is well defined, monotone, non-zero except at zero,
and itself satisfies \eqref{eq:admis} (Lemma~\ref{cor:U}).
Thus, $U_\mathcal{M}$ provides
the optimal distribution-free normalisation for the class:
\[
	\sqrt{\eta_f / U_{\mathcal{M}}(\eta_f)}, \qquad \text{$0<\eta_f<\infty$} \,.
\]
The main concern of this paper is a uniform and sharp bound on the optimal
adjustment function $\sqrt{U_\mathcal{M}}$. We prove the following result.

\begin{theorem}[Main theorem]\label{thm:main}
	For every $x > 1$,
	\[
		\left| \sqrt{U_{\mathcal{M}}(x)} - \sqrt{\log_4 x} \right|
		\le 2 \min\!\left\{ 1, \,
		\frac{2 \ln(\e + \ln x)}{\sqrt{\ln x}} \right\} \,.
	\]
\end{theorem}

A useful consequence is the following bound:
\[
	\mathbb{P}\left\{\frac{\sqrt{\eta_f}}{\sqrt{\log_4 \eta_f} +2}
	\sup_{t \in \mathbb{R}}
	\left| f_Z(t) - \Exf{f_Z(t)} \right|
	\ge \varepsilon \right\}
	\le 2\e^{-2\varepsilon^2}, \qquad \text{for all } \varepsilon > 0 \,,
\]
which holds for all $f_Z\in \mathcal{M}$ with finite $\eta_f\ge  1$.

Members of the class $\mathcal M$ are the normal forms of randomised signed-measure
processes that are uniformly norm bounded and order bounded from below,
in $\mathrm{ca}(\Omega,\Sigma)$, the Banach lattice of finite
countably additive signed measures. In this
representation, boundedness from below  is the structure that prevents
the pathologies relevant to the concentration problem
(Proposition~\ref{prop:maharam_cumulative} and Corollary~\ref{cor:maharam}).
The argument here isolates the
effect of this order boundedness on sharp sub-Gaussian normalisation
within the concentration-of-measure and empirical-process literature
\cite{boucheron2013concentration, ledoux2001concentration,
	talagrand1996new, vanderVaartWellner1996, vershynin2018high, wainwright2019high}.

As in Massart's proof, the upper bound proceeds by a uniformisation argument.
A multi-step parameter-side construction reduces the problem to a canonical subclass of
$\mathcal{M}$ without essentially reducing the required penalty below $U_{\mathcal{M}}$.
On this subclass an envelope of convex gauges exposes an elementary ramp-propagation property.
Tuning the envelope at each effective sample size gives a continuous upper bound for
the right-continuous version of $U_{\mathcal{M}}$, and this is consistent with the bound
of the main theorem.

The lower bound is finite-sample and combinatorial. To each subset of the
Hamming cube we associate a randomised monotone function whose deviation
is controlled by the covering radius and whose effective sample size by
the traversal length. Suitable sparse subsets give the matching logarithmic
bound of the main theorem.

Directly comparable are the standard bracketing bounds for the subclass of $\mathcal{M}$
consisting of averages of monotone summands~\cite{vanderVaartWellner1996}.
For some instances of the subclass, those bounds are exponentially loose
relative to Theorem~\ref{thm:main}.
The reason is that those bounds normalise by $\sqrt{n}$,
where $n$ is the number of independent summands,
rather than by $\sqrt{\eta_f}$. For the empirical distribution
with scalar data, the two
normalisations coincide, because $\eta_{F_n} = n$. However, the effective
sample size may far outstrip the number of independent summands: under the
strict stratification of Section~\ref{sec:EDFcluster}, where each of the
$n$ independent blocks contributes $m_n$ dependent ordered observations,
we have $\eta_f = n m_n^2$, and the bracketing normalisation is loose by a factor
that grows with the block size (Remark~\ref{rem:bracketing_loose}).
This same construction attains the sharpness of Theorem~\ref{thm:main}:
the configuration of averages of monotone summands on which bracketing is loose is the one on
which our logarithmic bound is sharp.

The smooth Gaussian process of Diebolt and Posse~\cite{diebolt1996density}
is another point of comparison. The suprema of such processes satisfy a
tail bound governed by the $\ell_2$ arc length of the trajectory on the
sphere, through the expected number of level crossings. The
distribution-free analogue, an order-bounded signed-measure process with
smooth non-monotone paths, is treated in
Corollary~\ref{cor:differentiable_paths}. The two bounds have the same
form, with the trajectory traced on the $\ell_1$-sphere rather than the
$\ell_2$-sphere; the $\ell_1$ and $\ell_2$ arc lengths differ in value.
The accompanying factor is a fixed constant in the Gaussian case, whereas
in the distribution-free setting it is the sharp logarithmic penalty of
the main theorem. Monotonicity buys Gaussian-like geometry across all
distributions, at the cost of boundedness and the logarithmic penalty.

The remainder of the paper is organised as follows.
Section~\ref{sec:equivalences} establishes the basic properties
of $U_{\mathcal{M}}$, including the equivalence framework and the
measurability of the centred supremum. Section~\ref{sec:lower_bound}
proves the lower bound using Hamming-cube embeddings, tuned admissible
constellations, and an interpolation estimate for the resulting
discrete sequence. Section~\ref{sec:uniformisation} develops a
parameter-side uniformisation, the classical data-side approach being
unavailable here. Section~\ref{sec:upper_bound} uses it to prove the
upper bound through a propagation property of the monotone functions
in $\mathsf{M}$. The proof of the symmetric bound in
Theorem~\ref{thm:main} is assembled in Section~\ref{sec:proof_main}.
Section~\ref{sec:applications} gives two applications. We first present
an empirical distribution function for stratified clustered data that
realises the lower-bound construction of the main theorem. We then
extend the bounds to non-monotone processes, bounding the suprema of
smooth bounded processes by their $\ell_1$ arc length. The Appendix
collects auxiliary results and deferred proofs.

\clearpage

\section{Preliminaries}\label{sec:equivalences}

This section establishes the following preliminary result.

\begin{lemma}\label{cor:U}
	The function $U_{\mathcal{M}}$ is well defined and monotone on
	$\mathbb{R}_+$; $U_{\mathcal M}(0)=0$, $U_{\mathcal M}(x)>0$ for all $x>0$, and
	\[
		U_{\mathcal{M}}(x) \ge  \frac{1}{2 \ln 2}, \qquad \text{for every } x \ge 1\,.
	\]
	Moreover, for every $f_Z \in \mathcal{M}$ with $\eta_f < \infty$, we have
	\[
		\mathbb{P}\left\{ \sqrt{\eta_f}\, \sup_{t \in \mathbb{R}}
		\left| f_Z(t) - \Exf{f_Z(t)} \right|
		\ge \varepsilon \sqrt{U_{\mathcal{M}}(\eta_f)} \right\}
		\le 2 \e^{-2\varepsilon^2} , \qquad \text{for all } \varepsilon > 0 \,.
	\]
\end{lemma}

The proof draws on several lemmas used throughout the paper. We end the
section by characterising $\mathcal{M}$ as the normal form of
signed-measure processes that are order bounded from below; the proof of this
characterisation uses the machinery developed in this section.

\subsection{General equivalence}

The following notion of equivalence is parameter-side, modulo the
expectation function.

\begin{definition}\label{def:equivalence}
	Two randomised functions $f_Z \colon \Theta \to \mathbb{R}$ and
	$g_Z \colon T \to \mathbb{R}$ over a common product measurable space
	$\mathcal{Z}$ and random vector $Z$ are \emph{equivalent} if
	$\eta_f=\eta_g=\infty$ or
	\[
		\sqrt{\eta_f}\, \sup_{\theta \in \Theta}
		\left| f(z, \theta) - \Exf{f_Z(\theta)} \right|
		=
		\sqrt{\eta_g}\, \sup_{t \in T}
		\left| g(z, t) - \Exf{g_{Z}(t)} \right|,
		\qquad \text{for each } z \in \mathcal{Z}\,.
	\]
\end{definition}

We give a construction that preserves equivalence,
and that will often be used.

\begin{lemma}\label{lem:elementary}
	Let $T$ and $\Theta^\ast$ be non-empty sets, and let $\Theta\subseteq\Theta^\ast$ be non-empty.
	Let $\psi\colon T\to\Theta^\ast$ be onto, let $\ds\colon T\to\mathbb R$
	be a bounded deterministic function, and let $a\neq0$.
	If $f_Z \colon \Theta \to \mathbb{R}$ is a randomised function, then
	\[
		g(z,t) :=
		\begin{cases}
			a\,f(z,\psi(t))+\ds(t), & \psi(t)\in\Theta\,,    \\
			\ds(t),                 & \psi(t)\notin\Theta\,,
		\end{cases}
	\]
	defines a randomised function $g_Z \colon T \to \mathbb{R}$ equivalent
	to $f_Z$; moreover, $a^2 \eta_g = \eta_f$.
\end{lemma}

\begin{proof}
	Deterministic terms cancel after centring, $\psi$ covers
	$\Theta$, and the coordinate coefficients satisfy $c_i(g)=|a|c_i(f)$.
	This implies $a^2 \eta_g = \eta_f$ and that the normalised centred supremum
	is unchanged.
\end{proof}

\subsection{General adjustments}

We now define the optimal adjustment $u_{f_Z}$ attached to a single
randomised function.

\begin{definition}\label{def:adjustment}
	For an arbitrary randomised function $f_Z \colon \Theta \to \mathbb{R}$ with finite
	$\eta_f$ whose centred supremum
	\[
		\sup_{\theta \in \Theta} \left| f_Z(\theta) - \Exf{f_Z(\theta)} \right|
	\]
	is measurable, denote by $u_{f_Z}$ the infimum of all $u \ge 0$ satisfying
	\begin{equation*}
		\mathbb{P}\left\{ \sqrt{\eta_f} \sup_{\theta \in \Theta}
		\left| f_Z(\theta) - \Exf{f_Z(\theta)} \right|
		\ge \varepsilon \sqrt{u} \right\}
		\le 2\e^{-2\varepsilon^2}, \qquad \text{for all } \varepsilon > 0 \,.
	\end{equation*}
	Any $u$ satisfying the above condition is called \emph{admissible} for $f_Z$.
\end{definition}

Equivalent randomisations in the sense of Definition~\ref{def:equivalence},
with finite $\eta$ and measurable centred supremum, have pointwise equal
normalised centred suprema, and thus the same adjustment.

Next we record the optimality property. It is immediate because bounded random
variables are sub-Gaussian.

\begin{lemma}\label{lem:opt1}
	Let $X$ be a bounded non-negative random variable, and let $A$ be the set of
	$u \ge 0$ satisfying
	\[
		\mathbb{P}\{ X \ge \varepsilon \sqrt{u} \}
		\le 2 \e^{-2\varepsilon^2}, \qquad \text{for all } \varepsilon > 0 \,.
	\]
	If $\mathbb{P}\{X > 0\} > 0$, then $A = [u^\ast, \infty)$, where
	\[
		u^\ast = \sup_{\left\{t > 0 \colon \mathbb{P}\{X \ge t\} > 0\right\}}
		\frac{2t^2}{\ln\left(\frac{2}{\mathbb{P}\{X \ge t\}}\right)}\,,
	\]
	which is positive, and otherwise $A = (0,\infty)$.
\end{lemma}

The next result is a consequence of the previous
characterisation of the admissible set.

\begin{lemma}\label{lem:u-upper}
	If $f_Z \colon \Theta \to [0,1]$ is a randomised function with finite
	$\eta_f$ and measurable centred supremum
	\[
		\xi =  \sup_{\theta \in \Theta} |f_Z(\theta) - \Exf{f_Z(\theta)}|\,,
	\]
	then $u_{f_Z} \le 2\eta_f / \ln 2$. Further,
	$\mathbb{P}\{\xi > 0\} > 0$ if and only if $u_{f_Z} > 0$, in
	which case $u_{f_Z}$ is admissible.
\end{lemma}

\subsection{Measurability of the centred supremum}

We show the measurability of the centred supremum of randomised monotone functions.

\begin{lemma}\label{lem:measurability}
	For every $f_Z \in \mathcal{M}$, the map
	\[
		z \mapsto \sup_{t \in \mathbb{R}}
		\left| f(z, t) - \Exf{f_Z(t)} \right|
	\]
	is measurable. In particular, $u_{f_Z}$
	is well defined for every $f_Z \in \mathcal{M}$ with finite $\eta_f$.
\end{lemma}

\begin{proof}
	We use the following elementary fact, whose proof is a
	standard exercise in monotonicity: if $g, h \in \mathsf{M}$ and
	$Q_g := \mathbb{Q} \cup D_g$, where $D_g$ is the countable set of
	discontinuities of $g$, then
	\[
		\sup_{t \in \mathbb{R}} \left| h(t) - g(t) \right|
		=
		\sup_{t \in Q_g} \left| h(t) - g(t) \right|\,.
	\]
	Let $f_Z \in \mathcal{M}$, and write $\mathbf{F}_Z(t) := \Exf{f_Z(t)}$
	for the expectation of the randomised monotone function. Because
	$f_z \in \mathsf{M}$ for each $z$, the function $\mathbf{F}_Z$ also
	belongs to $\mathsf{M}$, as a pointwise expectation of members of
	$\mathsf{M}$. Hence, for each $z \in \mathcal{Z}$,
	\[
		\sup_{t \in \mathbb{R}} \left| f_z(t) - \mathbf{F}_Z(t) \right|
		=
		\sup_{t \in Q_{\mathbf{F}_Z}} \left| f_z(t) - \mathbf{F}_Z(t) \right|\,.
	\]
	The set $Q_{\mathbf{F}_Z}$ is countable and does not depend on $z$.
	For each $t \in Q_{\mathbf{F}_Z}$, the map
	\[
		z \mapsto |f(z, t) - \mathbf{F}_Z(t)|
	\]
	is measurable, because
	$z \mapsto f(z, t)$ is measurable and $\mathbf{F}_Z(t)$ does not
	depend on $z$. Therefore,
	\[
		z \mapsto \sup_{t \in Q_{\mathbf{F}_Z}}
		\left| f(z, t) - \mathbf{F}_Z(t) \right|
	\]
	is measurable, being the supremum of a countable family of measurable
	functions, and the displayed equality completes the proof.
\end{proof}

\subsection{Richness of the monotone class}

The class $\mathcal M$ contains scaled indicators of coordinate maxima,
and these realise every effective sample size.

\begin{lemma}\label{lem:well-def}
	For every $x > 0$ there exists $f_Z \in \mathcal{M}$ with $\eta_f = x$ and $u_{f_Z}>0$.
\end{lemma}

\begin{proof}
	Choose $n \in \mathbb{N}$ with $nx \ge 1$, let $\mathcal{Z}_i = \{0, 1\}$,
	and let $Z_1, \ldots, Z_n$ be independent Bernoulli$(1/2)$ $\{0,1\}$-valued random
	variables. Define
	\[
		f_Z(t) := \frac{1}{\sqrt{nx}}
		\prod_{i=1}^n \mathbf{1}_{\{Z_i \le t\}}\,.
	\]
	Because $nx \ge 1$, the function $f_Z$ takes values in $[0,1]$ and is
	non-decreasing, so $f_Z \in \mathcal{M}$. Flipping a single coordinate
	$Z_i$ alters the product of indicators by at most $1$, with equality
	attained, so $c_i(f) = 1/\sqrt{nx}$. Summing the squares gives
	$\eta_f = \bigl( n \cdot (1/\sqrt{nx})^2 \bigr)^{-1} = x$.

	For every realisation $z$, the function $f_z$ is not equal to its expectation
	function, so the centred supremum is positive. By Lemma~\ref{lem:measurability}
	this supremum is measurable, and hence
	\[
		\mathbb P\left\{\sup_t |f_Z(t)-\Exf{f_Z(t)}|>0 \right\}=1\,.
	\]
	Lemma~\ref{lem:u-upper} then gives $u_{f_Z}>0$.
\end{proof}

Returning to Definition~\ref{def:adjustment}, we shall say that a function
$U\colon\mathbb R_+\to\mathbb R_+$ is \emph{pointwise admissible} if $U(\eta_f)$ is admissible
for every $f_Z\in\mathcal M$ with $\eta_f<\infty$.

\subsection{Admissibility for the monotone class}

We restate the optimality property for the monotone class.

\begin{lemma}\label{lem:opt}
	For any $f_Z \in \mathcal{M}$ with finite $\eta_f$, every $u \ge u_{f_Z}$
	with $u > 0$ is admissible for $f_Z$.
\end{lemma}

\begin{proof}
	By Lemma~\ref{lem:measurability} the centred supremum is measurable, so
	Lemma~\ref{lem:u-upper} applies. If $u_{f_Z} = 0$, the centred supremum
	is almost surely zero, and every $u > 0$ is admissible.
	If \mbox{$u_{f_Z} > 0$}, then \mbox{$u_{f_Z}$} is admissible by
	Lemma~\ref{lem:u-upper} and so is $u$.
\end{proof}

\begin{proof}[Proof of Lemma~\ref{cor:U}]
	Define $U \colon \mathbb{R}_+ \to \mathbb{R}_+$ by $U(0) = 0$ and, for
	$x > 0$,
	\[
		U(x) := \sup\{u_{f_Z} : f_Z \in \mathcal{M},\ \eta_f = x\}\,.
	\]
	By Lemma~\ref{lem:well-def} the set is non-empty, and
	by Lemma~\ref{lem:u-upper} it is bounded above by $2x/\ln 2$,
	so $U$ is well defined.

	\medskip
	\emph{Admissibility of $U$.}  The value at zero is irrelevant for admissibility,
	because no $f_Z\in\mathcal M$ has $\eta_f=0$. For $x>0$,
	$U(x) \ge u_{f_Z}$ for every $f_Z \in \mathcal{M}$ with $\eta_f=x$.
	Because $x>0$, there exists by Lemma~\ref{lem:well-def}
	$f_Z\in \mathcal{M}$ with $\eta_f=x$ and $u_{f_Z}>0$. Thus, $U(x)>0$ and is
	admissible by Lemma~\ref{lem:opt}.

	\medskip
	\emph{Identification $U_\mathcal{M} = U$.} For any pointwise admissible
	$V \colon \mathbb{R}_+ \to \mathbb{R}_+$ and any $f_Z \in \mathcal{M}$
	with $\eta_f = x$, the bound at $V(x)$ holds, so $V(x) \ge u_{f_Z}$.
	Taking the supremum over $f_Z$ gives $V(x) \ge U(x)$, and taking the
	infimum over pointwise admissible $V$ gives $U_\mathcal{M}(x) \ge U(x)$.
	Conversely, $U$ is itself pointwise admissible by the previous paragraph, so
	$U_\mathcal{M}(x) \le U(x)$. Hence $U_\mathcal{M} = U$, and the
	displayed inequality of the lemma is the admissibility
	of $U_\mathcal{M}$ at $\eta_f$.

	\medskip
	\emph{Monotonicity.} Let $y \ge x > 0$ and let $f_Z \in \mathcal{M}$
	with $\eta_f = x$. Define $g_Z := \sqrt{x/y}\, f_Z$. Because $0<\sqrt{x/y}\le 1$,
	the function $g_Z$ still takes values in $[0,1]$ and is non-decreasing.
	Thus, $g_Z \in \mathcal{M}$. By Lemma~\ref{lem:elementary} with $a = \sqrt{x/y}$,
	$\eta_g = y$ and $g_Z$ is equivalent to $f_Z$, hence
	$u_{g_Z} = u_{f_Z}$. Taking the supremum over $f_Z$ with $\eta_f = x$
	gives $U_\mathcal{M}(y) \ge U_\mathcal{M}(x)$.

	\medskip
	\emph{Lower bound.} Specialise the construction of
	Lemma~\ref{lem:well-def} to $n = 1$, so $f_Z(t) = \mathbf{1}_{\{Z_1 \le t\}}$
	with $\eta_f = 1$, and let $Z_1$ be uniform on $\{0, 1\}$. For each
	$t \in [0, 1)$, $f_Z(t) = \mathbf{1}_{\{Z_1 = 0\}}$ takes the values $0$
	and $1$ with equal probability, so
	$|f_Z(t) - \Exf{f_Z(t)}| = \tfrac{1}{2}$. Hence
	$\sqrt{\eta_f}\, \sup_t |f_Z(t) - \Exf{f_Z(t)}| = \tfrac{1}{2}$
	surely.

	Lemma~\ref{lem:opt1} gives
	$u_{f_Z} = 1/(2\ln 2)$. Hence
	\[
		U_\mathcal{M}(1) \ge 1/(2 \ln 2)\,,
	\]
	and monotonicity extends this to every $x \ge 1$.
\end{proof}

\subsection{Characterisation of the class by randomised signed-measure processes} \label{sec:cummu}

We describe a two-way representation of the class $\mathcal{M}$ in terms of
signed-measure processes that are order bounded from below and indexed by
regular chains.

Let $(\Omega,\Sigma)$ be a measurable space. Consider an
indexed family $\{E_\theta\}_{\theta\in\Theta}\subseteq\Sigma$ such that
\begin{enumerate}
	\item for all $\theta,\theta^\ast\in\Theta$, either
	      \[
		      E_\theta\subseteq E_{\theta^\ast}
		      \qquad\text{or}\qquad
		      E_{\theta^\ast}\subseteq E_\theta\,;
	      \]
	\item there exists a positive, finite, finitely additive measure $\lambda$ on
	      $(\Omega,\Sigma)$ such that
	      \[
		      \theta\ne\theta^\ast
		      \quad\Longrightarrow\quad
		      \lambda(E_\theta\triangle E_{\theta^\ast})>0\,.
	      \]
\end{enumerate}
We call such an indexed family a \emph{regular chain}.

Condition~(2) states that there exists $\lambda^*$ that is positive on the algebra
generated by the increments of the chain, its restriction to the increments is $\lambda$.
Kelley~\cite{kelley1959} characterises the existence of such a
$\lambda^*$ combinatorially, through intersection numbers, and equivalently
through strictly positive linear functionals on the associated
Banach lattice; $\lambda$ plays the latter role here.

Let $\{\mu_z\}_{z\in\mathcal Z}$ be a family of
finite countably additive signed measures in
the Dedekind complete Banach lattice $\mathrm{ca}(\Omega,\Sigma)$
endowed with the total variation norm.
Assume that the family is order bounded from below and uniformly norm bounded,
\[
	\|\mu_z\|\le 1,\qquad z\in\mathcal Z\,.
\]
Assume further that, for each fixed $\theta\in\Theta$, the map
\[
	z\longmapsto \mu_z(E_\theta)
\]
is measurable.

For a random vector $Z$ taking values in $\mathcal Z$ with independent
components, define the randomised function
\[
	F_Z\colon\Theta\to\mathbb R, \qquad F_Z(\theta):=\mu_Z(E_\theta)\,.
\]
We call $F_Z$ the \emph{cumulative function} of the randomised
signed measure $\mu_Z$ along the regular chain $\{E_\theta\}_{\theta\in\Theta}$.

The proof of the following result is in the Appendix.

\begin{proposition}\label{prop:maharam_cumulative}
	Under the hypotheses above, the cumulative function $F_Z$ is
	equivalent to a member $f_Z$ of $\mathcal M$. The member $f_Z$ may be
	chosen with effective sample size
	\[
		\eta_f
		= \bigl(1+\|\nu\|\bigr)^2\,\eta_F,
		\qquad \nu:=\bigwedge_{z\in\mathcal Z}\mu_z\,.
	\]
	Conversely, every member of $\mathcal{M}$ is equivalent to a
	cumulative function of a randomised signed measure, along a regular chain,
	satisfying the hypotheses above with $\nu=0$.
\end{proposition}

The basic idea is that $\mathbf{1}_{E_\theta}$, understood as linear functionals,
form a chain of distinct points in the order dual of $\mathrm{ca}(\Omega,\Sigma)$,
and $\lambda$ is a linear functional that is strictly positive on the
convex cone generated by the increments of this chain.
The converse is proved by a standard argument using the lexicographic
ordering of $\mathbb{R}^3$.

Together with the one-sided upper bound (Theorem~\ref{thm:upper_bound_main}) and Lemma~\ref{cor:U},
the proposition gives our inequality, in explicit form, for randomised signed-measure processes.

\begin{corollary}\label{cor:maharam}
	For every cumulative function $F_Z$ as above with $1<\eta_F<\infty$,
	\[
		\mathbb{P}\left\{\frac{\sqrt{\eta_F}}
		{\mathcal{R}\bigl((1+\|\nu\|)^2\,\eta_F\bigr)}
		\sup_{\theta \in \Theta}
		\left| F_Z(\theta) - \Exf{F_Z(\theta)} \right|
		\ge \varepsilon \right\}
		\le 2\e^{-2\varepsilon^2}, \qquad \text{for all } \varepsilon > 0 \,,
	\]
	where $\nu:=\bigwedge_{z\in\mathcal Z}\mu_z$ and
	\[
		\mathcal{R}(x) := \sqrt{\log_4 x} +
		2\min\!\left\{1,\, \frac{\ln(\e+\ln x)}{\sqrt{\ln x}}\right\},
		\qquad x>1 \,.
	\]
\end{corollary}

\begin{proof}
	By Proposition~\ref{prop:maharam_cumulative}, choose an equivalent
	$f_Z\in\mathcal M$ with
	\[
		\eta_f=(1+\|\nu\|)^2\eta_F\,.
	\]
	Equivalence gives equality of the normalised centred suprema. Because
	$\eta_f>1$, Theorem~\ref{thm:upper_bound_main} and
	Lemma~\ref{cor:U} imply
	\[
		\sqrt{U_{\mathcal M}(\eta_f)}
		\le \mathcal R(\eta_f)\,,
	\]
	which gives the displayed bound.
\end{proof}

\section{The lower bound}\label{sec:lower_bound}

We prove the following theorem.

\begin{theorem}\label{thm:lower_bound_main}
	For all $x > 1$,
	\[
		\left( \sqrt{U_{\mathcal{M}}(x)} - \sqrt{\log_4 x} \right)^- \le
		2 \min\left\{1, \frac{2\ln(\e+ \ln x)}{\sqrt{\ln x}}\right\} \,.
	\]
\end{theorem}

The proof uses Lemma~\ref{lem:penalty_bound}, which applies to every
non-empty subset of the Hamming cube and converts its covering radius and
shortest-path traversal length into a lower bound on
$\sqrt{U_{\mathcal{M}}(x)}$.

\subsection{The Hamming cube}

For $n \in \mathbb{N}$, let $\mathcal{Q}_n := \{0,1\}^n$ denote the discrete
cube equipped with the Hamming distance
$d_H(q,q') = \#\{i \colon q_i \neq q'_i\}$. For any non-empty
\mbox{$S \subseteq \mathcal{Q}_n$}, define the covering radius
\[
	\rho(S) = \max_{q \in \mathcal{Q}_n} \min_{s \in S} d_H(q, s) \,,
\]
and the shortest-path traversal length
\[
	\mathcal{L}(S) = \min \left\{ \sum_{j=2}^{|S|} d_H(s_j, s_{j-1}) \colon
	(s_1, s_2, \dots, s_{|S|}) \text{ is an ordering of } S \right\} \,.
\]

The following result ties Hamming-cube geometry to randomised monotone
functions.

\begin{lemma}\label{lem:penalty_bound}
	For every non-empty $S \subseteq \mathcal{Q}_n$,
	\[
		\sqrt{U_{\mathcal{M}}\left( \frac{(\mathcal{L}(S) + n)^2}{n} \right)}
		\ge \frac{n - 2\rho(S)}{\sqrt{2n \ln 2}} \,.
	\]
\end{lemma}

We record a consequence.

\begin{corollary}\label{cor:naive_bound}
	For all $x \ge 4$,
	\[
		\sqrt{U_{\mathcal{M}}(x)} \ge \frac{1}{\sqrt{2\ln 2}} \sqrt{\log_4 x - 1} \,.
	\]
\end{corollary}

\begin{proof}
	Let $x \ge 4$ and define $n_x := \lfloor \log_4 x \rfloor$. Set
	$S = \mathcal{Q}_{n_x}$. We have $\rho(S) = 0$. The hypercube $\mathcal{Q}_{n_x}$
	admits a Hamiltonian path with $\mathcal{L}(S) \le  2^{n_x} - 1$.

	We have
	\[
		(2^{n_x} - 1 + n_x)^2 \le n_x 4^{n_x}\,.
	\]
	Indeed,
	\[
		2^{n_x}+n_x-1 \le \sqrt{n_x}\,2^{n_x}\,.
	\]
	For $n_x=1$ this is equality. For $n_x\ge2$, it follows because
	\[
		n_x-1=(\sqrt{n_x}-1)(\sqrt{n_x}+1)
		\le (\sqrt{n_x}-1)2^{n_x}\,.
	\]
	Thus, because $\mathcal{L}(S)\le 2^{n_x}-1$ and $4^{n_x}\le x$,
	\[
		\frac{(\mathcal{L}(S)+n_x)^2}{n_x} \le 4^{n_x} \le x \,.
	\]

	Lemma~\ref{lem:penalty_bound} and the monotonicity of $U_{\mathcal{M}}$ imply
	\[
		\sqrt{U_{\mathcal{M}}(x)} \ge \frac{n_x}{\sqrt{2n_x\ln 2}} =
		\frac{1}{\sqrt{2\ln 2}}\sqrt{n_x} \,.
	\]
	Because $n_x \ge \log_4 x - 1$, the claim follows.
\end{proof}

\begin{proof}[Proof of Lemma~\ref{lem:penalty_bound}]
	If $n - 2\rho(S) \le 0$, the lower bound is non-positive and holds
	because $U_{\mathcal{M}}$ is non-negative. Assume $n - 2\rho(S) > 0$.

	Identify $\{0,1\}^n$ with the oriented cube $\mathcal{Z} = \{-1,1\}^n$ via
	the standard bijection. Identify $S$ with a subset $\Theta \subseteq
		\mathcal{Z}$, preserving $\rho(\Theta) = \rho(S)$ and $\mathcal{L}(\Theta) =
		\mathcal{L}(S)$. Define
	\begin{equation}\label{eq:f-extreme}
		f(z,\theta) = \langle z, \theta \rangle \,.
	\end{equation}

	Because $\langle z, \theta \rangle = n - 2d_H(z,\theta)$, we have
	\[
		\min_{z\in \mathcal{Z}} \max_{\theta\in \Theta} \langle z, \theta \rangle
		= n - 2\rho(\Theta) \,.
	\]
	This establishes $\sup_{\theta\in \Theta} f(z,\theta) \ge n - 2\rho(\Theta)$
	for all $z\in \mathcal{Z}$.

	We have $c_i(f) = 2$ for all $i$, because changing $z_i$ to $-z_i$
	alters $\langle z, \theta\rangle$ by $\pm 2\theta_i$ with
	\mbox{$\theta_i \in \{-1, 1\}$}. This gives $\eta_f = 1/(4n)$. Thus,
	\[
		\sqrt{\eta_f} \sup_{\theta\in \Theta} f(z,\theta)
		\ge \frac{n - 2\rho(\Theta)}{2\sqrt{n}}, \qquad \text{for every } z \in \mathcal{Z} \,.
	\]

	Let $Z$ be uniformly distributed over $\mathcal{Z}$. Because
	$\Exf{f_Z(\theta)}=0$ and $n - 2\rho(\Theta) > 0$, we have
	\[
		\sqrt{\eta_f} \sup_{\theta\in\Theta} |f(z,\theta) -
		\Exf{f_Z(\theta)}| \ge \sqrt{\eta_f} \max_{\theta\in\Theta}
		f(z,\theta) \ge \frac{n - 2\rho(\Theta)}{2\sqrt{n}}
	\]
	for every $z\in\mathcal{Z}$. In particular,
	\[
		\mathbb{P}\left\{ \sqrt{\eta_f} \sup_{\theta\in \Theta}
		\left| f_Z(\theta) - \Exf{f_Z(\theta)} \right|
		\ge \frac{n - 2\rho(\Theta)}{2\sqrt{n}} \right\} = 1 \,.
	\]
	Therefore,
	\[
		\sqrt{u_{f_Z}} \ge \frac{n - 2\rho(\Theta)}{\sqrt{2n \ln 2}} \,.
	\]
	To see this, note that for arbitrary $0<u < (n - 2\rho(\Theta))^2/(2n\ln 2)$,
	taking
	\[
		\varepsilon =  (n - 2\rho(\Theta))/(2\sqrt{n}\sqrt{u})
	\]
	gives $\varepsilon\sqrt{u} = (n - 2\rho(\Theta))/(2\sqrt{n})$ and
	$2\e^{-2\varepsilon^2} < 1$. Thus, $u_{f_Z}\ge u$.

	Let $m = |\Theta|$ and order $\Theta = \{\theta_1, \theta_2, \ldots, \theta_m\}$
	to match the shortest path $\mathcal{L}(\Theta)$. Define the surjective
	mapping $\sigma \colon \mathbb{R} \to \Theta$
	by $\sigma(t) = \theta_{j(t)}$, where
	\[
		j(t) = \max\{1, \min\{m, \lfloor t \rfloor\}\} \,.
	\]
	Define $g \colon \mathcal{Z} \times \mathbb{R} \to \mathbb{R}$ by
	\[
		g(z,t) = \frac{1}{2(\mathcal{L}(\Theta) + n)}\left( n +
		f(z, \sigma(t)) + 2\sum_{k=2}^{j(t)} d_H(\theta_k, \theta_{k-1}) \right) \,,
	\]
	where the sum evaluates to zero for $j(t)=1$.

	Because $f(z, \theta) \in [-n, n]$ and the sum of successive distances is at
	most $\mathcal{L}(\Theta)$, the mapping $g_Z$ takes values in $[0,1]$.
	If $j\ge2$, then
	\[
		f(z,\theta_j)-f(z,\theta_{j-1})
		\ge -2d_H(\theta_j,\theta_{j-1})\,.
	\]
	The added term, therefore, compensates for any fall in the inner product, so
	$g_Z(t)$ is non-decreasing in $t$. Thus, $g_Z \in \mathcal{M}$.

	By Lemma~\ref{lem:elementary}, $g_Z$ is equivalent
	to $f_Z$ with
	\[
		\eta_{g} = 4(\mathcal{L}(\Theta) + n)^2 \eta_f
		= \frac{(\mathcal{L}(\Theta)+n)^2}{n} \,.
	\]
	Evaluating $U_{\mathcal{M}}$ gives
	\[
		\sqrt{U_{\mathcal{M}}\left(\frac{(\mathcal{L}(\Theta)+n)^2}{n}\right)}
		= \sqrt{U_{\mathcal{M}}(\eta_g)} \ge \sqrt{u_{g_Z}} = \sqrt{u_{f_Z}}
		\ge \frac{n - 2\rho(\Theta)}{\sqrt{2n \ln 2}} \,,
	\]
	completing the proof.
\end{proof}

\subsection{The deterministic range}

We prove a consequence of Lemma~\ref{cor:U} and
Corollary~\ref{cor:naive_bound}.

\begin{lemma}\label{lem:deterministic_regime}
	For all $1 < x \le 4^{111}$, the bound in Theorem~\ref{thm:lower_bound_main}
	is satisfied.
\end{lemma}

\begin{proof}
	For $1 < x \le 4^8$, the bound holds by the lower estimate in
	Lemma~\ref{cor:U}. Specifically, we have
	$U_{\mathcal{M}}(x) \ge (2\ln 2)^{-1}$ for all $x\ge 1$. Note that
	\[
		\sqrt{\log_4(4^8)} - (2\ln 2)^{-1/2} = \sqrt{8} - (2\ln 2)^{-1/2} < 2\,,
	\]
	and
	\[
		2 \min\left\{1, \frac{2\ln(\e+ \ln x)}{\sqrt{\ln x}}\right\} = 2,
		\qquad \text{for } 1< x\le 4^{55}\,.
	\]
	Together these establish the bound for $1 < x \le 4^8$.

	By Corollary~\ref{cor:naive_bound}, for $x \ge 4^8$, the residual gap is
	bounded by
	\[
		\sqrt{\log_4 x} - \sqrt{U_{\mathcal{M}}(x)}  \le
		\sqrt{\log_4 x} - \frac{1}{\sqrt{2\ln 2}}\sqrt{\log_4 x - 1} \,.
	\]
	Letting $w = \log_4 x$ and $a = 2\ln 2$, direct calculation shows that
	\[
		\sqrt{w} - \frac{\sqrt{w-1}}{\sqrt{a}} \le 2, \qquad
		1 \le w \le 170\,,
	\]
	and
	\[
		\sqrt{a}w -\sqrt{w^2-w} \le 4\ln(\e+aw), \qquad
		1 \le w \le 111\,.
	\]
	These together give
	\[
		\sqrt{w} - \frac{\sqrt{w-1}}{\sqrt{a}}
		\le 2\min\left\{1, \frac{2\ln(\e+aw)}{\sqrt{aw}}\right\} \,,
	\]
	for the remaining range $8 \le w \le 111$.
\end{proof}

\subsection{The admissible constellations}

Admissible constellations are the integer parameter systems used to extend
the lower bound beyond the deterministic range.

\begin{definition}
	A tuple $(k,t_k,n_k,r_k,m_k)\in\mathbb{N}^5$ is a
	\emph{constellation} if
	\[
		\begin{gathered}
			n_k = k\,t_k^2, \qquad
			r_k = \frac{k\,t_k(t_k-1)}{2}, \\[2mm]
			t_k \ge 2\,.
		\end{gathered}
	\]
	It is an \emph{admissible constellation} if $m_k\le 2^{n_k}$.
\end{definition}

Here, $k$ is the sequence index, $t_k$ is the \emph{tuning} parameter, $n_k$ is
the dimension of the Hamming cube, $r_k$ is the target covering radius, and
$m_k$ is the cardinality of the target subset within the $n_k$-dimensional cube.
The associated \emph{volume} is $n_k m_k^2$.

The next lemma records that admissible constellations are generated by the
three integer parameters $k$, $t_k$, and $m_k$.

\begin{lemma}\label{lem:constellation_from_parameters}
	Let $k,t_k\in\mathbb{N}$ with $t_k\ge2$. Define
	\[
		n_k := k\,t_k^2,
		\qquad
		r_k := \frac{k\,t_k(t_k-1)}{2}\,.
	\]
	Hence $n_k,r_k\in\mathbb{N}$. Consequently, for every
	$m_k\in\mathbb{N}$ satisfying $m_k\le 2^{n_k}$, the tuple
	\[
		(k,t_k,n_k,r_k,m_k)
	\]
	is an admissible constellation.
\end{lemma}

\begin{proof}
	That $n_k$ is an integer is immediate. The number $r_k$ is an integer because
	$t_k(t_k-1)$ is even. The remaining defining conditions are $t_k\ge2$ and
	$m_k\le2^{n_k}$.
\end{proof}

We now apply Lemma~\ref{lem:penalty_bound} to admissible
constellations.

\begin{lemma}\label{lem:summarise}
	For an admissible constellation $(k, t_k, n_k, r_k, m_k)$,
	if there exists a non-empty $S \subseteq \mathcal{Q}_{n_k}$ satisfying
	$|S| \le m_k$ and $\rho(S) \le r_k$, then
	\[
		\sqrt{U_{\mathcal{M}}(n_k m_k^2 )} \ge \frac{1}{\sqrt{2\ln 2}} \sqrt{k} \,.
	\]
\end{lemma}

\begin{proof}
	The shortest-path length is bounded by
	\[
		\mathcal{L}(S) \le n_k(|S|-1) \,,
	\]
	which implies $(\mathcal{L}(S) + n_k)^2/n_k \le n_k|S|^2$. Thus, the
	monotonicity of $U_{\mathcal{M}}$ and Lemma~\ref{lem:penalty_bound} imply
	\begin{align*}
		\sqrt{U_{\mathcal{M}}(n_k m_k^2 )}
		 & \ge \sqrt{U_{\mathcal{M}} (n_k |S|^2 )}                                        \\
		 & \ge \sqrt{U_{\mathcal{M}} \left( \frac{(\mathcal{L}(S) + n_k)^2}{n_k} \right)} \\
		 & \ge \frac{1}{\sqrt{2\ln 2}} \frac{n_k - 2\rho(S)}{\sqrt{n_k}}                  \\
		 & \ge \frac{1}{\sqrt{2\ln 2}} \frac{n_k - 2r_k}{\sqrt{n_k}} \,.
	\end{align*}
	Because any constellation satisfies
	\[
		\frac{n_k - 2r_k}{\sqrt{n_k}} = \sqrt{k} \,,
	\]
	the stated bound is established.
\end{proof}

\subsection{Probabilistic existence}

A standard probabilistic argument establishes sufficient conditions for the
existence of the subset required in Lemma~\ref{lem:summarise}.

\begin{lemma}\label{lem:existence_condition_weakened}
	For an admissible constellation $(k, t_k, n_k, r_k, m_k)$,
	if
	\[
		\sum_{j=0}^{r_k} \binom{n_k}{j} \ge \frac{2^{n_k}}{m_k} n_k \ln 2 \,,
	\]
	then
	\[
		\sqrt{U_{\mathcal{M}} (n_k m_k^2 )} \ge \frac{1}{\sqrt{2\ln 2}} \sqrt{k} \,.
	\]
\end{lemma}

\begin{proof}
	Let $s_1, \dots, s_{m_k}$ be independent random variables uniformly distributed
	over $\mathcal{Q}_{n_k}$. Applying the union bound over all
	$q \in \mathcal{Q}_{n_k}$, we obtain
	\begin{align*}
		\mathbb{P}\left\{ \max_{q\in \mathcal{Q}_{n_k}} \min_{1 \le j \le m_k}
		d_H(q,s_j) > r_k \right\}
		 & \le 2^{n_k}\left(1 - \frac{1}{2^{n_k}}
		\sum_{j=0}^{r_k}\binom{n_k}{j}\right)^{m_k}  \\
		 & < 2^{n_k} \exp\left[- \frac{m_k}{2^{n_k}}
			\sum_{j=0}^{r_k}\binom{n_k}{j}\right] \,.
	\end{align*}
	Under the stated hypothesis, this probability is strictly less than $1$. This
	establishes the existence of a set $S \subseteq \mathcal{Q}_{n_k}$ satisfying
	$|S| \le m_k$ and $\rho(S) \le r_k$. The conclusion follows directly from
	Lemma~\ref{lem:summarise}.
\end{proof}

Retaining only one binomial term, at the cost of slack, gives a closed-form
sufficient condition for the preceding lemma.

\begin{lemma}\label{lem:m}
	For an admissible constellation $(k, t_k, n_k, r_k, m_k)$,
	if
	\[
		\ln(m_k^2) \ge k + 3\ln k +
		\left[6\ln t_k + \frac{k}{6(t_k^2-1)}\right]
		+ \left[\ln 2 + 2\ln(\ln 2)\right] \,,
	\]
	then
	\[
		\sqrt{U_{\mathcal{M}}(n_k m_k^2)} \ge \frac{1}{\sqrt{2\ln 2}} \sqrt{k} \,.
	\]
\end{lemma}

\begin{proof}
	By Lemma~\ref{lem:existence_condition_weakened}, the conclusion holds if
	\[
		\frac{1}{2^{n_k}} \sum_{j=0}^{r_k} \binom{n_k}{j} \ge \frac{1}{m_k} n_k \ln 2 \,.
	\]
	Lemma~\ref{lem:true_binomial_term_general} in
	Appendix~\ref{sec:appendix_binomial}
	provides a lower bound on the relevant binomial shell:
	\[
		\frac{1}{2^{n_k}}\binom{n_k}{r_k} \ge \frac{1}{\sqrt{2k}\,t_k}
		\exp\left(-\frac{k}{12(t_k^2-1)}\right) \e^{-k/2} \,.
	\]
	Bounding the summation from below by this single term, it suffices that
	\[
		\frac{1}{\sqrt{2k}\,t_k}
		\exp\left(-\frac{k}{12(t_k^2-1)}\right) \e^{-k/2}
		\ge \frac{1}{m_k} n_k \ln 2 \,.
	\]
	Substituting $n_k = k t_k^2$, taking logarithms, and isolating $\ln m_k$,
	we obtain
	\begin{align*}
		\ln m_k
		 & \ge \ln(k t_k^2 \ln 2) + \frac{1}{2}\ln(2k) + \ln t_k +
		\frac{k}{12(t_k^2-1)} + \frac{k}{2}                                      \\
		 & = \frac{k}{2} + \frac{3}{2}\ln k + 3\ln t_k + \frac{k}{12(t_k^2-1)} +
		\frac{1}{2}\ln 2 + \ln(\ln 2) \,.
	\end{align*}
	Multiplying by $2$ establishes the stated bound.
\end{proof}

\subsection{The tuned constellations}\label{subsec:calibrated}

The sufficient condition in Lemma~\ref{lem:m} is equivalently a lower
bound on $m_k$. We define the real parameter $\tilde m_k$ by taking
equality in its logarithmic form:
\[
	2\ln(\tilde{m}_k) = k + 3\ln k +
	\left[6\ln t_k + \frac{k}{6(t_k^2-1)}\right]
	+ \left[\ln 2 + 2\ln(\ln 2)\right] \,.
\]
Applying the identity $n_k = k t_k^2$, we obtain
\begin{align*}
	\ln(n_k \tilde{m}_k^2)
	 & = k + 4\ln k + 8\ln t_k + \frac{k}{6(t_k^2-1)}
	+ \left[\ln 2 + 2\ln(\ln 2)\right]                                 \\
	 & < k + 4\ln k + \left[8\ln t_k + \frac{k}{6(t_k^2-1)}\right] \,,
\end{align*}
where the strict inequality follows from $\ln 2 + 2\ln(\ln 2) < 0$.

Minimising the dominant terms of the tuning penalty
\[
	8\ln t_k + \frac{k}{6(t_k^2-1)}
\]
gives
\[
	t_k \sim \frac{\sqrt{k}}{\sqrt{24}} \,.
\]
We set the calculation-friendly real tuning parameter
\[
	\tilde{t}_k = \frac{\sqrt{k}}{5} \,.
\]

\begin{definition}
	A constellation $(k, t_k, n_k, r_k, m_k)$ is \emph{tuned}
	if
	\[
		t_k = \left\lceil \tilde{t}_k \right\rceil
		\quad\text{and}\quad
		m_k = \lceil \tilde{m}_k \rceil \,.
	\]
	If it is also admissible, then it is an \emph{admissible tuned constellation}.
\end{definition}

We now establish the existence of admissible tuned constellations.

\begin{lemma}\label{lem:admissibility}
	For all $k \ge 26$, there exists a unique admissible tuned constellation
	\[
		(k, t_k, n_k, r_k, m_k) \,.
	\]
	Furthermore, this constellation satisfies
	\[
		\sqrt{U_{\mathcal{M}}(n_k m_k^2 )} \ge \frac{1}{\sqrt{2\ln 2}} \sqrt{k} \,.
	\]
\end{lemma}

\begin{proof}
	Fix $k \ge 26$. The formulas
	\[
		t_k=\left\lceil\tilde t_k\right\rceil,
		\qquad
		m_k=\left\lceil\tilde m_k\right\rceil
	\]
	uniquely determine $n_k$ and $r_k$. It remains to verify admissibility.

	Because $\tilde t_k=\sqrt{k}/5$ and $k\ge26$, we have
	\[
		2\le t_k\le\sqrt{k}\,.
	\]
	Note also that $\tilde m_k>1$. By its defining equation and $t_k\ge2$, we have
	\[
		2\ln\tilde m_k
		\ge k+3\ln k+\ln2+2\ln(\ln2)>0
	\]
	for $k\ge26$.

	We have
	\begin{align*}
		\ln m_k
		 & = \ln \lceil \tilde{m}_k \rceil \le \ln(\tilde{m}_k + 1)
		\le \ln(2\tilde{m}_k) = \ln \tilde{m}_k + \ln 2                        \\
		 & < \frac{k}{2} + \frac{3}{2}\ln k + 3\ln t_k + \frac{k}{12(t_k^2-1)}
		+ \ln 2                                                                \\
		 & \le \frac{k}{2} + \frac{3}{2}\ln k + 3\ln t_k + \frac{k}{36}
		+ \ln 2                                                                \\
		 & \le \frac{k}{2} + \frac{3}{2}\ln k + 3\ln \sqrt{k} + \frac{k}{36}
		+ \ln 2                                                                \\
		 & = \frac{19}{36}k + 3\ln k + \ln 2                                   \\
		 & \le \frac{19}{36}k + k                                              \\
		 & = \frac{55}{36}k < 4k\ln 2 \le k t_k^2 \ln 2 = n_k\ln 2\,.
	\end{align*}
	This establishes $m_k \le 2^{n_k}$, and hence the tuned constellation is
	admissible.

	Also, $m_k=\lceil\tilde m_k\rceil\ge\tilde m_k$, so
	$\ln(m_k^2)\ge\ln(\tilde m_k^2)$. By the definition of $\tilde m_k$,
	the hypothesis of Lemma~\ref{lem:m} is satisfied, and the final statement
	follows from that lemma.
\end{proof}

\subsection{The probabilistic range}

For each $k \ge 26$, let
\[
	(k, t_k, n_k, r_k, m_k)
\]
be the unique admissible tuned constellation, the existence of which is
guaranteed by Lemma~\ref{lem:admissibility}.

The next result locates each $x>4^{111}$ between two consecutive tuned volumes.

\begin{lemma}\label{lem:wedge}
	For all $x > 4^{111}$, there exists an integer $k_x \ge 120$ satisfying
	\[
		n_{k_x} m_{k_x}^2 \le x < n_{k_x+1} m_{k_x+1}^2 \,.
	\]
\end{lemma}

\begin{proof}
	Fix $x > 4^{111}$ and define the index set
	\[
		\mathcal{I}_x = \{k \ge 26 \colon n_k m_k^2 \le x \} \,.
	\]
	Because the sequence $(n_k m_k^2)$ diverges to infinity,
	the set $\mathcal{I}_x$ has finite cardinality. We now show
	that it is not empty.

	At the index $k = 120$, we have $t_{120} = 3$, and $n_{120} = 1080$.
	Expanding the logarithm gives
	\begin{align*}
		\ln(n_{120} m_{120}^2)
		 & \le \ln(n_{120}(\tilde{m}_{120}+1)^2)        \\
		 & = \ln(n_{120}\tilde{m}_{120}^2) +
		2\ln\left(1 + \frac{1}{\tilde{m}_{120}}\right)  \\
		 & < 150.4 + 0.6 \le 151 < 222\ln 2 < \ln x \,.
	\end{align*}
	Thus, $120 \in \mathcal{I}_x$, so the set is non-empty.

	Define $k_x = \max \mathcal{I}_x$. The maximality of $k_x \ge 120$
	ensures $k_x+1 \notin \mathcal{I}_x$, establishing the
	bounds
	\[
		n_{k_x} m_{k_x}^2 \le x < n_{k_x+1} m_{k_x+1}^2 \,,
	\]
	and completing the proof.
\end{proof}

\subsection{Estimating the gap}

Combining Lemmas~\ref{lem:admissibility} and \ref{lem:wedge}, for
each $x > 4^{111}$ we have
\begin{align*}
	\sqrt{\log_4 x} - \sqrt{U_{\mathcal{M}}(x)}
	 & \le \frac{1}{\sqrt{2\ln 2}} \left( \sqrt{\ln x} - \sqrt{k_x} \right) \\
	 & = \frac{1}{\sqrt{2\ln 2}} \left( \frac{\ln x - k_x}
	{\sqrt{\ln x} + \sqrt{\ln x - (\ln x - k_x)}} \right) \,,
\end{align*}
where $k_x \ge 120$ is the index specified in Lemma~\ref{lem:wedge}.
The rearrangement is valid because $k_x\le \ln x$. Indeed,
$m_{k_x}\ge\tilde m_{k_x}$, and the definition of $\tilde m_{k_x}$ gives
\[
	\ln(n_{k_x}m_{k_x}^2)
	\ge \ln(n_{k_x}\tilde m_{k_x}^2)
	\ge k_x\,.
\]
Because $n_{k_x}m_{k_x}^2\le x$, it follows that $k_x\le\ln x$.

We now estimate the gap $\ln x - k_x$ using an interpolation based on
Lemma~\ref{lem:wedge}.

\begin{lemma}\label{lem:laststep}
	For $x > 4^{111}$, we have
	\[
		\ln x - k_x < 8\ln(\ln x) \,.
	\]
\end{lemma}

\begin{proof}
	Let $y = \ln x$. By Lemma~\ref{lem:wedge}, we have
	\begin{align*}
		y
		 & \ge \ln(n_{k_x} m_{k_x}^2) \ge \ln(n_{k_x} \tilde{m}_{k_x}^2)  \\
		 & = k_x + 4 \ln k_x + 8 \ln t_{k_x} + \frac{k_x}{6(t_{k_x}^2-1)}
		+ \ln 2 + 2\ln(\ln 2)                                             \\
		 & > k_x + 4 \ln k_x + 8 \ln t_{k_x} + \ln 2 + 2\ln(\ln 2) \,.
	\end{align*}
	Substituting the lower bounds $k_x \ge 120$ and $t_{k_x} \ge 3$ establishes
	$y > k_x + 27$. This verifies the strict inequality
	\begin{equation}\label{eq:1kx}
		y > k_x + 1 \,.
	\end{equation}

	Define $K_x = k_x + 1$. By Lemma~\ref{lem:wedge}, we have
	\begin{align*}
		y
		 & \le \ln(n_{K_x} m_{K_x}^2)
		\le \ln(n_{K_x} \tilde{m}_{K_x}^2) + \frac{2}{\tilde{m}_{K_x}}      \\
		 & < K_x + 4 \ln K_x + 8 \ln t_{K_x} + \frac{K_x}{6(t_{K_x}^2-1)} +
		\frac{2}{\tilde{m}_{K_x}} \,.
	\end{align*}
	We bound the two tuning penalties. By definition,
	\[
		t_{K_x} \le \frac{\sqrt{K_x}}{5} + 1 =
		\frac{\sqrt{K_x}}{5}\left(1 + \frac{5}{\sqrt{K_x}}\right) \,.
	\]
	Combining this with the elementary inequality $\ln(1+z) \le z$ for $z > -1$
	gives
	\[
		8 \ln t_{K_x} \le 4 \ln K_x - 8 \ln 5 + \frac{40}{\sqrt{K_x}} \,.
	\]
	Furthermore, the lower bound $t_{K_x} \ge \frac{\sqrt{K_x}}{5}$ ensures
	\[
		\frac{K_x}{6(t_{K_x}^2-1)} \le \frac{25 K_x}{6(K_x-25)} =
		\frac{25}{6}\left(1 + \frac{25}{K_x-25}\right) \,.
	\]

	For $K_x\ge121$, the first four terms in the bracket below are bounded from below by
	their values at $K_x=121$. Also, $\tilde m_{K_x}>1$, so
	$2/\tilde m_{K_x}<2$. Direct calculation gives
	\[
		8\ln5-1-\frac{40}{11}
		-\frac{25}{6}\left(1+\frac{25}{96}\right)
		-2
		> 0 \,.
	\]
	Substituting $K_x = k_x + 1$ and evaluating the constants for $K_x \ge 121$
	gives
	\begin{align*}
		y - k_x
		 & \le 1 + 8 \ln K_x + \frac{40}{\sqrt{K_x}} - 8\ln 5 +
		\frac{25}{6}\left(1 + \frac{25}{K_x-25}\right) + \frac{2}{\tilde{m}_{K_x}} \\
		 & \le 8 \ln(k_x+1) - \left[ 8\ln 5 - 1 - \frac{40}{\sqrt{K_x}} -
			\frac{25}{6}\left(1 + \frac{25}{K_x-25}\right) - \frac{2}{\tilde{m}_{K_x}}
		\right]                                                                    \\
		 & < 8\ln(k_x+1) \,.
	\end{align*}
	Applying the logarithmic substitution $\ln(k_x+1) < \ln y = \ln(\ln x)$ derived
	from \eqref{eq:1kx} establishes the stated result.
\end{proof}

\begin{proof}[Proof of Theorem~\ref{thm:lower_bound_main}]
	For $1 < x \le 4^{111}$ the bound holds by
	Lemma~\ref{lem:deterministic_regime}. We show the bound for
	$x > 4^{111}$.

	Fix $x > 4^{111}$ and let $y = \ln x$. By
	Lemmas~\ref{lem:admissibility} and~\ref{lem:wedge},
	\[
		\sqrt{\log_4 x} - \sqrt{U_{\mathcal{M}}(x)} \le
		\frac{1}{\sqrt{2\ln 2}} \left( \frac{y - k_x}
		{\sqrt{y} + \sqrt{y - (y - k_x)}} \right) \,.
	\]
	The map $u \mapsto \frac{u}{\sqrt{y} + \sqrt{y - u}}$ is strictly
	increasing for $0 \le u < y$. By Lemma~\ref{lem:laststep} and the
	condition $y > 222\ln 2$,
	\[
		y - k_x < 8\ln y < y \,,
	\]
	so applying the monotonicity above,
	\begin{align*}
		\frac{1}{\sqrt{2\ln 2}} \frac{y - k_x}{\sqrt{y} + \sqrt{y - (y - k_x)}}
		 & < \frac{1}{\sqrt{2\ln 2}} \frac{8\ln y}{\sqrt{y} + \sqrt{y - 8\ln y}} \\
		 & = \frac{8\ln y}{\sqrt{y}} \left[ \sqrt{2\ln 2}
			\left( 1 + \sqrt{1 - \frac{8\ln y}{y}} \right) \right]^{-1} \,.
	\end{align*}
	The function $y \mapsto \frac{\ln y}{y}$ is strictly decreasing for
	$y \ge \e$. For $y > 222\ln 2$, the bracketed factor in the
	inverted expression exceeds $2$. Therefore,
	\[
		\frac{1}{\sqrt{2\ln 2}} \frac{y - k_x}{\sqrt{y} + \sqrt{y - (y - k_x)}}
		< \frac{4\ln y}{\sqrt{y}} \,.
	\]
	Because $\ln y < \ln(\e + y) = \ln(\e + \ln x)$, this
	establishes the bound stated in
	Theorem~\ref{thm:lower_bound_main} for $x > 4^{111}$:
	for $y>222\ln2$, we also have $2\ln(\e+y)\le\sqrt y$,
	so the minimum in the theorem equals its second argument.
\end{proof}

\section{Parameter-side uniformisation}\label{sec:uniformisation}

Massart's~\cite{massart1990tight} upper-bound argument proceeds by
uniformising the data. That approach is not available for a general
randomised monotone function. The construction in this section is reminiscent of
Debreu's gap lemma~\cite{debreu1964continuity}: the expectation function
is made strictly increasing, and its gaps are filled by
interpolation. Equivalence is preserved throughout.

Let $\mathcal{M}_{\mathsf{id}}([0,1])$ be the set of randomised monotone
functions $f_Z\in\mathcal M$ satisfying the section-wise endpoint conditions
\begin{equation}\label{eq:edge}
	f_z(t)=0 \quad (t\le 0),\qquad f_z(t)=1 \quad (t\ge 1)\,,
\end{equation}
and the diagonal expectation condition
\[
	\Exf{f_Z(t)} = t,\qquad t\in[0,1] \,.
\]
Because the functions are deterministic outside $[0,1]$, we shall often
identify randomised functions satisfying \eqref{eq:edge} with their
restrictions to $[0,1]$. Conversely, any randomised function on $[0,1]$
whose endpoints are deterministically equal to $0$ and $1$, respectively, may
be identified with its constant extension satisfying \eqref{eq:edge}.

Let $U_{\mathcal{M}_{\mathsf{id}}}(0)=0$, and for each $x>0$ define
\[
	U_{\mathcal{M}_{\mathsf{id}}}(x) := \sup\left\{ u_{f_Z} \colon
	f_Z \in \mathcal{M}_{\mathsf{id}}([0,1]), \ \eta_f = x \right\} \,.
\]
The properties of $U_{\mathcal{M}_{\mathsf{id}}}$ are established in the
theorem of the present section.

\begin{theorem}\label{thm:unif}
	The function $U_{\mathcal{M}_{\mathsf{id}}}\colon \mathbb{R}_+\to \mathbb{R}_+$
	is well defined and monotone.  Further, the right-continuous versions of
	$U_{\mathcal{M}_{\mathsf{id}}}$ and $U_{\mathcal{M}}$ are equal everywhere.
\end{theorem}

The proof requires an interpolation argument.

\subsection{Extension by partition of unity}

Let $\Lambda$ be an arbitrary non-empty set and
\[
	f \colon \mathcal{Z}\times \Theta \times \Lambda \to [0,1]
\]
be such that $\theta \mapsto f_Z(\theta\mid\lambda)$ is a randomised
function for each fixed $\lambda \in \Lambda$. Here $\,\cdot\mid\lambda$
denotes selection of the $\lambda$-section at the fixed, non-random
$\lambda\in\Lambda$.

Assume that $\Lambda$ has a majorant $\lambda^\ast \in \Lambda$ satisfying
\[
	c_i(f(\cdot, \cdot\mid \lambda^\ast)) \ge c_i(f(\cdot, \cdot\mid\lambda))\,,
\]
for all $i \in \{1, \dots, n\}$ and all $\lambda \in \Lambda$, and
\[
	\sup_{\theta\in \Theta} \left|f(z,\theta\mid\lambda^\ast)-
	\Exf{f_Z(\theta\mid\lambda^\ast)}\right| \ge
	\sup_{\theta\in \Theta} \left|f(z,\theta\mid\lambda)- \Exf{f_Z(\theta\mid\lambda)}\right|\,,
\]
for all $z\in \mathcal{Z}$ and all $\lambda\in\Lambda$.

Let $T$ be a set satisfying $\Theta \subseteq T$, and let
$\beta \colon T \to [0,1]^{\Theta\times \Lambda}$ be a partition of unity
satisfying:
\begin{enumerate}
	\item $\beta_{(\theta, \lambda)}(t) = 0$, for all but finitely many
	      $(\theta,\lambda)$, and
	      $\sum_{(\theta,\lambda)\in \Theta\times \Lambda} \beta_{(\theta,\lambda)} (t)
		      = 1$ for all $t \in T$.
	\item For every $\theta \in \Theta$ we have
	      $\beta_{(\theta, \lambda^\ast)}(\theta) = 1$.
\end{enumerate}

Define the function $g_Z \colon T \to [0,1]$ by
\[
	g_Z(t) := \sum_{(\theta,\lambda) \in \Theta\times \Lambda}
	f_Z(\theta\mid\lambda)\, \beta_{(\theta, \lambda)}(t) \,.
\]

\begin{lemma}\label{lem:discrete_interpolation}
	The function $g_Z$ is a randomised function equivalent to
	$f_Z(\cdot\mid\lambda^\ast)$, and
	\[
		c_i(g) = c_i(f(\cdot,  \cdot\mid \lambda^\ast)) \,,
	\]
	for all $i \in \{1, \dots, n\}$.
\end{lemma}

\begin{proof}
	Because $\beta_{(\cdot,\cdot)}(\cdot)$ is a partition of unity, $g_Z(t)$ is a finite convex
	combination of the variables $f_Z(\theta\mid\lambda) \in [0,1]$; in particular,
	$z \mapsto g_Z(t)$ is measurable as a finite sum of measurable sections, so
	$g_Z$ is a randomised function.

	For any coordinate $i$, the triangle inequality ensures the variation of the sum is
	bounded by the convex combination of the variations $c_i(f(\cdot,\cdot\mid\lambda))$,
	each at most $c_i(f(\cdot, \cdot\mid\lambda^\ast))$, giving
	\[
		c_i(g) \le c_i(f(\cdot, \cdot\mid \lambda^\ast))\,.
	\]

	Similarly, centring the function and applying the triangle inequality gives,
	for each $z \in \mathcal{Z}$,
	\[
		\sup_{t \in T} \left|g(z,t) - \Exf{g_Z(t)}\right| \le
		\sup_{\theta \in \Theta} \left|f(z,\theta\mid\lambda^\ast) -
		\Exf{f_Z(\theta\mid\lambda^\ast)}\right|\,.
	\]

	Equality in both inequalities follows from the second property of the
	partition of unity and the non-negativity of the weights.
	For any $\theta \in \Theta$,
	evaluating at $t = \theta$ gives $\beta_{(\theta, \lambda^\ast)}(\theta) = 1$ and
	$\beta_{(\theta', \lambda')}(\theta) = 0$ for all other pairs
	$(\theta', \lambda')$. Thus, $g_Z(\theta) = f_Z(\theta\mid\lambda^\ast)$.

	So $c_i(g) \ge c_i(f(\cdot,\cdot\mid\lambda^\ast))$, and, because the supremum over
	$T$ is bounded below by the supremum over the subset $\Theta \subseteq T$,
	exact equality holds.
\end{proof}

\subsection{The uniformisation}

The theorem of the section is an immediate consequence of the following
lemma, whose proof is constructive.

\begin{lemma}\label{lem:uniformisation}
	For any randomised function $f_Z$ in $\mathcal{M}$ and any $a>0$, there is an
	equivalent randomised function $g_Z$ in $\mathcal{M}_{\mathsf{id}}([0,1])$
	satisfying
	\[
		(1+a)^{-2}\eta_g = \eta_f \,.
	\]
\end{lemma}

\begin{proof}

	\medskip
	\noindent\emph{Restricting the domain.}
	Fix $f_Z \in \mathcal{M}$ and let $a > 0$. Let $\psi \colon (1/4,3/4) \to
		\mathbb{R}$ be a strictly increasing surjective function. Define
	\[
		f^{(1)}(z,t) :=
		\begin{cases}
			0,                       & t \le 1/4   \,,  \\
			f\left(z,\psi(t)\right), & t\in(1/4,3/4)\,, \\
			1,                       & t \ge 3/4 \,.
		\end{cases}
	\]
	By Lemma~\ref{lem:elementary}, and because outside $(1/4,3/4)$ the function is
	constant, the randomised function $f^{(1)}_Z\in \mathcal{M}$ is equivalent
	to $f_Z$.

	\medskip
	\noindent\emph{Strictly increasing expectations.}
	Restricting $f^{(1)}$ to the domain $[0,1]$, define the affine scaling
	\[
		f^{(2)}(z,t) := \frac{1}{1+a} f^{(1)}(z,t) + \frac{a}{1+a} t \,,\qquad
		t\in[0,1] \,.
	\]
	By Lemma~\ref{lem:elementary}, the function $f^{(2)}_Z$ is equivalent to
	$f^{(1)}_Z$, and
	\[
		(1+a)^{-2} \eta_{f^{(2)}} = \eta_{f^{(1)}} \,.
	\]
	Hence, $\eta_{f^{(2)}} = (1+a)^2 \eta_f$. Because $a>0$, the added linear drift
	ensures the mapping \mbox{$t \mapsto f^{(2)}(z,t)$} is strictly increasing for
	every $z$. Consequently, the expectation function $\Exf{f^{(2)}_Z(t)}$ is strictly
	increasing.

	\medskip
	\noindent\emph{Generalised inverse.}
	For every $\vartheta\in [0,1]$, define the inverse mapping
	\[
		t_\vartheta := \inf \left\{ t \in [0,1] \colon \vartheta \le
		\Exf{f^{(2)}_Z(t)} \right\} =
		\sup \left\{ t \in [0,1] \colon \Exf{f^{(2)}_Z(t)} \le \vartheta \right\} \,.
	\]
	Let $\Theta \subseteq [0,1]$ be the set of values $\theta$ satisfying
	$\theta = \Exf{f^{(2)}_Z(t_\theta)}$.

	\medskip

	\noindent\emph{Target interpolation.}
	If the expectation function is discontinuous, then $\Theta$ has gaps. We use
	interpolation to fill these gaps and extend the domain from $\Theta$ to $[0,1]$.

	Define the left- and right-continuous versions $f_Z^{(2)\mathrm{L}}$ and
	$f_Z^{(2)\mathrm{R}}$, respectively, by
	\[
		f^{(2)\mathrm{L}}(z,x) := \lim_{y\uparrow x} f^{(2)}(z,y)
		\quad \text{and}\quad
		f^{(2)\mathrm{R}}(z,x) := \lim_{y\downarrow x} f^{(2)}(z,y)\,,
	\]
	for each $(z,x)$. By Lemma~\ref{lem:RL}, $f_Z^{(2)\mathrm{L}}$ and
	$f_Z^{(2)\mathrm{R}}$ are monotone and randomised functions satisfying
	\[
		c_i(f^{(2)\mathrm{L}}) = c_i(f^{(2)\mathrm{R}})\le c_i(f^{(2)})\,,
	\]
	for all $i$. Further,
	\begin{align*}
		\sup_{x\in [0,1]}\left| f^{(2)\mathrm{L}}(z,x)-\Exf{f^{(2)\mathrm{L}}_Z(x)}\right|
		 & =
		\sup_{x\in [0,1]}\left| f^{(2)\mathrm{R}}(z,x)-\Exf{f^{(2)\mathrm{R}}_Z(x)}\right| \\
		 & \le
		\sup_{x\in [0,1]}\left| f^{(2)}(z,x)-\Exf{f_Z^{(2)}(x)}\right|\,,
	\end{align*}
	for all $z$.

	Exactly one of the following must hold for any $\vartheta \in [0,1]$:
	\begin{align}
		\Exf{f^{(2)}_Z(t_\vartheta)}           & < \vartheta \le
		\Exf{f^{(2)\mathrm{R}}_Z(t_\vartheta)}\,, \tag{R}               \\
		\Exf{f^{(2)\mathrm{L}}_Z(t_\vartheta)} & \le \vartheta <
		\Exf{f^{(2)}_Z(t_\vartheta)}\,, \tag{L}                         \\
		\Exf{f^{(2)}_Z(t_\vartheta)}           & = \vartheta\,. \tag{E}
	\end{align}
	Thus, there is a unique weight vector
	$(\alpha_{\vartheta}, \alpha_{\mathrm{R} \vartheta},
		\alpha_{\mathrm{L} \vartheta}) \in \mathbb{R}^3_+$ satisfying the simplex
	conditions:
	\[
		\alpha_{\vartheta}+ \alpha_{\mathrm{R} \vartheta}+
		\alpha_{\mathrm{L} \vartheta}=1,
		\qquad
		\alpha_{\mathrm{L} \vartheta} \alpha_{\mathrm{R} \vartheta}=0,
		\qquad
		\alpha_{\mathrm{L} \vartheta} + \alpha_{\mathrm{R} \vartheta}=0
		\quad \text{in case } (E)\,,
	\]
	that interpolates the target expectation:
	\[
		\vartheta = \alpha_{\vartheta} \Exf{f^{(2)}_Z(t_\vartheta)} +
		\alpha_{\mathrm{R} \vartheta} \Exf{f^{(2)\mathrm{R}}_Z(t_\vartheta)} +
		\alpha_{\mathrm{L} \vartheta} \Exf{f^{(2)\mathrm{L}}_Z(t_\vartheta)} \,.
	\]
	Clearly $\alpha_\vartheta=1$ if and only if
	$\vartheta = \Exf{f^{(2)}_Z(t_\vartheta)}$. Also, in case (R),
	$\alpha_{\mathrm{L} \vartheta}=0$; in case (L), $\alpha_{\mathrm{R} \vartheta}=0$.

	\medskip

	Define the target randomised function $g_Z$ in
	$\mathcal{M}_{\mathsf{id}}([0,1])$ via
	\[
		g(z,\vartheta) := \alpha_{\vartheta} f^{(2)}(z,t_\vartheta)
		+ \alpha_{\mathrm{R} \vartheta} f^{(2)\mathrm{R}}(z,t_\vartheta)
		+ \alpha_{\mathrm{L} \vartheta} f^{(2)\mathrm{L}}(z,t_\vartheta)\,.
	\]
	By Lemma~\ref{lem:interpolation_monotonicity} of the Appendix, the map
	$\vartheta \mapsto g(z, \vartheta)$ is non-decreasing for each $z \in \mathcal{Z}$.
	Linearity of expectation implies $\Exf{g_Z(\vartheta)} = \vartheta$ for all
	$\vartheta \in [0,1]$, establishing $g_Z \in \mathcal{M}_{\mathsf{id}}([0,1])$.

	\medskip
	\noindent\emph{Application of the partition of unity.}
	Set $T=[0,1]\supseteq\Theta$, and define on $\Theta$ the three functions
	\[
		f^{(3)}_Z(\theta):=f^{(2)}_Z(t_\theta),\qquad
		f^{(3)\mathrm{R}}_Z(\theta):=f^{(2)\mathrm{R}}_Z(t_\theta),\qquad
		f^{(3)\mathrm{L}}_Z(\theta):=f^{(2)\mathrm{L}}_Z(t_\theta)\,.
	\]
	Because $\theta\mapsto t_\theta$ maps $\Theta$ onto $[0,1]$, $f^{(3)}_Z$ is
	equivalent to $f^{(2)}_Z$, $f^{(3)\mathrm{R}}_Z$ to $f^{(2)\mathrm{R}}_Z$, and
	$f^{(3)\mathrm{L}}_Z$ to $f^{(2)\mathrm{L}}_Z$; each reparametrised pair of equivalent functions
	has the same coefficients $c_i(\cdot)$ and centred supremum.

	Let $\Lambda=\{\mathrm{E},\mathrm{R},\mathrm{L}\}$, each label corresponding to
	the superscript-$(3)$ function carrying it. By Lemma~\ref{lem:RL} applied
	to the superscript-$(2)$ functions, and the preceding equivalences, $f^{(3)}$
	dominates $f^{(3)\mathrm{R}}$ and $f^{(3)\mathrm{L}}$ in every coefficient
	$c_i(\cdot)$ and in centred supremum, so
	$\lambda^\ast=\mathrm{E}$ is a majorant of $\Lambda$.

	Set $\theta_\vartheta:=\Exf{f^{(2)}_Z(t_\vartheta)}\in\Theta$. Define the partition
	of unity $\beta\colon[0,1]\to[0,1]^{\Theta\times\Lambda}$ to be
	supported on $\{\theta_\vartheta\}\times\Lambda$, with ordered triplet
	\[
		\big(\beta_{(\theta_\vartheta,\mathrm{E})}(\vartheta),\,
		\beta_{(\theta_\vartheta,\mathrm{R})}(\vartheta),\,
		\beta_{(\theta_\vartheta,\mathrm{L})}(\vartheta)\big)
		:= (\alpha_{\vartheta},\,\alpha_{\mathrm{R}\vartheta},\,\alpha_{\mathrm{L}\vartheta})\,,
	\]
	and zero on all other pairs. For each $\vartheta$ the three weights
	$\alpha_{\vartheta},\alpha_{\mathrm{R}\vartheta},\alpha_{\mathrm{L}\vartheta}$
	sum to one and are carried by the pairs $(\theta_\vartheta,\lambda)$,
	$\lambda\in\Lambda$, so $\beta$ has finite support and weights summing to one; moreover, for
	$\theta\in\Theta$ we have $\theta_\theta=\theta$, so case~(E) holds with
	$\alpha_\theta=1$, giving $\beta_{(\theta,\mathrm{E})}(\theta)=1$. Thus, $\beta$ is a
	partition of unity in the sense of Lemma~\ref{lem:discrete_interpolation}.

	Because $t_{\theta_\vartheta}=t_\vartheta$,
	\[
		g(z,\vartheta)=\alpha_{\vartheta} f^{(3)}(z,\theta_\vartheta)
		+\alpha_{\mathrm{R}\vartheta} f^{(3)\mathrm{R}}(z,\theta_\vartheta)
		+\alpha_{\mathrm{L}\vartheta} f^{(3)\mathrm{L}}(z,\theta_\vartheta)\,.
	\]
	By Lemma~\ref{lem:discrete_interpolation}, $g_Z$ is equivalent to
	$f^{(3)}_Z$, with $c_i(g)=c_i(f^{(3)})=c_i(f^{(2)})$. Equivalence is transitive;
	$g_Z$ is equivalent to $f_Z$.
\end{proof}

\begin{proof}[Proof of Theorem~\ref{thm:unif}]
	By Lemma~\ref{lem:uniformisation}, $U_{\mathcal{M}_{\mathsf{id}}}$
	is well defined on the whole of $\mathbb{R}_+$, because $U_\mathcal{M}$
	is well defined. Further,
	\[
		U_{\mathcal{M}_{\mathsf{id}}}(x) \le U_{\mathcal{M}}(x)
		\le U_{\mathcal{M}_{\mathsf{id}}}\big((1+a)^2 x\big)\,,
	\]
	for any $x\ge 0$ and any $a>0$. Because $a$ is arbitrary,
	$U_{\mathcal{M}_{\mathsf{id}}}$ is monotone and its right-continuous version
	equals the right-continuous version of $U_{\mathcal{M}}$.
\end{proof}

\section{The upper bound}\label{sec:upper_bound}

This section establishes the following theorem.

\begin{theorem}\label{thm:upper_bound_main}
	For every $x> 1$,
	\[
		\sqrt{U_{\mathcal{M}}(x)} - \sqrt{\log_4 x} \le
		2\min\left\{1,\frac{\ln(\e+\ln x)}{\sqrt{\ln x}}\right\} \,.
	\]
\end{theorem}

Because $U_{\mathcal{M}}$ is monotone, its right-continuous version
dominates it, and by Theorem~\ref{thm:unif} it, in turn, equals
the right-continuous version of $U_{\mathcal{M}_{\mathsf{id}}}$. As
the envelope on the right of the theorem is continuous, it
suffices to establish the bound for $U_{\mathcal{M}_{\mathsf{id}}}$.
We begin with a geometric property of deterministic non-decreasing functions.

\subsection{Ramp propagation of non-decreasing functions}

The deviation of a non-decreasing function $f\colon[0,1]\to[0,1]$ from the
identity satisfies the following property, whose straightforward
proof is deferred to the Appendix.

\begin{lemma}[Ramp property]\label{lem:ramp-sup}
	Let $f\colon[0,1]\to[0,1]$ be a non-decreasing function.

	If $\sup_{t\in[0,1]} (f(t)-t)^+ \ge \varepsilon$, then there exists
	$t^{*}\in[0,1]$ such that for all $0<\alpha\le \varepsilon$,
	\[
		(f(t^{*}+\alpha)-(t^{*}+\alpha))^+ \ge \varepsilon-\alpha \,.
	\]

	If $\sup_{t\in[0,1]} (f(t)-t)^- \ge \varepsilon$, then there exists
	$t^{*}\in[0,1]$ such that for all $0<\alpha\le \varepsilon$,
	\[
		(f(t^{*}-\alpha)-(t^{*}-\alpha))^- \ge \varepsilon-\alpha \,.
	\]
\end{lemma}

A pointwise deviation induces a forward-propagating ramp of deviations,
bounding a triangular region adjacent to the identity (Figure~\ref{fig:ramp}).
\begin{figure}[htbp]
	\centering
	\begin{tikzpicture}[scale=1.2]
		\draw[->, thick, gray] (-0.2,0) -- (5.5,0) node[right, text=black] {$t$};
		\draw[->, thick, gray] (0,-0.2) -- (0,5.5) node[above, text=black] {};
		\draw[thick, domain=0:5, smooth, blue] plot(\x,{1.1*\x+sin(\x r)+0.41})
		node[right, text=black] {};

		\draw[dashed, thick] (2,2) -- (2,3.519);
		\filldraw[fill=blue!10, draw=blue, thick] (2,2) -- (2,3.519) --
		(3.519,3.519) -- cycle;

		\filldraw (2,3.519) circle (2pt) node[above left, text=black] {$f(t^*)$};

		\draw[darkgray] (0.9,0) arc (0:45:0.9);
		\node[darkgray] at (1.25,0.27) {$45^\circ$};
		\draw[blue, thick] (2.25,3.519) -- (2.25,3.269) -- (2,3.269);
		\draw[thin, black] (0,0) -- (5,5) node[above right, text=black]{};
		\draw[blue, thick] (1.93,2.76) -- (2.07,2.76);          \draw[blue, thick] (2.76,3.449) -- (2.76,3.589);
	\end{tikzpicture}
	\caption{A forward-propagating ramp.}
	\label{fig:ramp}
\end{figure}

\subsection{Isolating the ramp property with gauges}

Let $\mathcal{H}$ be the cone of $C^1$, strictly convex, strictly increasing
Young functions $H\colon\mathbb{R}_+\to\mathbb{R}_+$ satisfying $H(0)=0$. For
$H\in\mathcal{H}$, let \mbox{$h:=H'$} denote its strictly increasing gauge.

\begin{corollary}
	\label{cor:ramp-deterministic}
	If $f\colon[0,1]\to[0,1]$ is a deterministic non-decreasing function and
	$r>0$, then
	\[
		r \sup_{t\in[0,1]} |f(t)-t| \;\le\; \inf_{H\in\mathcal{H}}
		H^{-1}\!\left(r\int_0^1 h\!\left(r\,|f(t)-t|\right)\,\dd t\right) \,.
	\]
\end{corollary}
\begin{proof}
	Assume $m := r \sup_{t\in[0,1]} |f(t)-t| > 0$; the case $m=0$ is trivial.
	By definition either $\sup_t (f(t)-t)^+ = m/r$ or $\sup_t (f(t)-t)^- = m/r$;
	consider the former, the latter being symmetric. By
	Lemma~\ref{lem:ramp-sup} there is $t^\ast$ with $t^\ast+\alpha\in[0,1]$ for
	$0\le\alpha\le m/r$ and
	\[
		r\bigl(f(t)-t\bigr)\;\ge\; m-r\,(t-t^\ast),
		\qquad t\in (t^\ast,\,t^\ast+m/r] \,.
	\]
	Fix $H\in\mathcal{H}$ with gauge $h=H'\ge 0$ strictly increasing.
	Restricting the integral to this interval and using monotonicity of $h$,
	\begin{align*}
		r\int_0^1 h\!\left(r|f(t)-t|\right)\dd t
		 & \ge r\int_{t^\ast}^{t^\ast+m/r}
		h\!\bigl(m-r(t-t^\ast)\bigr)\dd t    \\
		 & = \int_0^m h(s)\,\dd s = H(m) \,,
	\end{align*}
	by the substitution $s=m-r(t-t^\ast)$. Because $H^{-1}$ is increasing,
	\[
		H^{-1}\!\left(r\int_0^1 h\!\left(r|f(t)-t|\right)\dd t\right)\ge m \,.
	\]
	As $H\in\mathcal{H}$ was arbitrary, the infimum is at least $m$.
\end{proof}

\subsubsection{Pointwise consequence for randomised functions}

We record the following corollary, whose proof is now immediate.

\begin{corollary}
	\label{cor:ramp-deterministic2}
	If $f_Z\in\mathcal{M}_{\mathsf{id}}([0,1])$ and $0<\eta_f<\infty$, then for
	every realisation $z$ of $Z$, we have
	\[
		\sqrt{\eta_f}\,\sup_{t\in[0,1]} |f_z(t)-t| \;\le\; \inf_{H\in\mathcal{H}}
		H^{-1}\!\left(\sqrt{\eta_f}\int_0^1
		h\!\left(\sqrt{\eta_f}\,|f_z(t)-t|\right)\dd t\right) \,.
	\]
\end{corollary}

\begin{proof}
	Fix a realisation $z$.
	The function $f_z\colon[0,1]\to[0,1]$ is deterministic and non-decreasing,
	so Corollary~\ref{cor:ramp-deterministic} applies with $r=\sqrt{\eta_f}$ and $f=f_z$.
\end{proof}

\subsection{The expected supremum in
\texorpdfstring{$\mathcal{M}_{\mathsf{id}}$}{M\_id}}

For any $f_Z \in \mathcal{M}_{\mathsf{id}}([0,1])$, $\Exf{f_Z(t)} = t$.
By McDiarmid's inequality, the pointwise absolute deviation satisfies the
sub-Gaussian tail:
\[
	\sup_{t\in[0,1]} \mathbb{P}\left\{ \sqrt{\eta_f}|f_Z(t)-t| \ge
	\varepsilon \right\} \le 2\e^{-2\varepsilon^2}, \qquad
	\varepsilon > 0 \,.
\]
Let $\mathcal{H}_0 \subset \mathcal{H}$ denote the subclass of gauges
satisfying the growth condition
\[
	\lim_{\varepsilon \to \infty} h(\varepsilon)\e^{-2\varepsilon^2} = 0 \,.
\]
Corollary~\ref{cor:ramp-deterministic2} provides a bound for the expected absolute
supremum.

\begin{lemma}\label{lem:random-supremum}
	For any $f_Z \in \mathcal{M}_{\mathsf{id}}([0,1])$ with finite $\eta_f$, we have
	\[
		\Exf{ \sqrt{\eta_f} \sup_{t\in[0,1]} |f_Z(t)-t| } \le \inf_{H\in\mathcal{H}_0}
		H^{-1}\left( \sqrt{\eta_f} \int_0^\infty 8\varepsilon\, h(\varepsilon)
		\e^{-2\varepsilon^2} \,\dd\varepsilon \right) \,.
	\]
\end{lemma}

\begin{proof}
	Set
	\[
		r:=\sqrt{\eta_f}, \qquad X_t(z):=r|f_z(t)-t|, \qquad X_t^R(z):=r|f_z^R(t)-t|,
	\]
	where $f_z^R$ is the right-continuous version of $f_z$ from Lemma~\ref{lem:RL}.

	By Lemma~\ref{lem:measurability}, $z\mapsto\sup_{t\in[0,1]}X_t(z)$
	and $z\mapsto\sup_{t\in[0,1]}X^R_t(z)$ are measurable.
	We write $X_t$ for $X_t(Z)$ and $X_t^R$ for $X_t^R(Z)$.

	Fix $H\in\mathcal H_0$, and write $h=H'$. By Corollary~\ref{cor:ramp-deterministic2}, for every $z\in \mathcal{Z}$,
	\[
		\sup_{t\in[0,1]}X_t(z) \le H^{-1}\left( r\int_0^1 h(X_t(z))\,\dd t \right)\,.
	\]
	For fixed $z$, the monotone functions $f_z$ and $f_z^R$ differ at most at the
	jump points of $f_z$, hence at most on a countable set. Therefore, we have the equality
	\[
		\int_0^1 h(X_t(z))\,\dd t = \int_0^1 h(X_t^R(z))\,\dd t \,,
	\]
	for every $z\in \mathcal{Z}$.  Additionally, $(z,t)\mapsto h(X_t^R(z))$ is jointly measurable.

	Taking expectations and then taking the infimum over $H\in\mathcal H_0$ gives
	\begin{align*}
		\Exf{\sup_{t\in[0,1]}X_t} & \le \inf_{H\in\mathcal H_0} \Exf{ H^{-1}\left( r\int_0^1 h(X_t)\,\dd t \right) }     \\
		                          & =  \inf_{H\in\mathcal H_0} \Exf{ H^{-1}\left( r\int_0^1 h(X^R_t)\,\dd t \right) }\,.
	\end{align*}
	Because $H^{-1}$ is concave and increasing, Jensen's inequality gives
	\[
		\Exf{ H^{-1}\left( r\int_0^1 h(X_t^R)\,\dd t \right) } \le H^{-1}\left( r\,\Exf{\int_0^1 h(X_t^R)\,\dd t} \right)\,.
	\]
	By joint measurability and non-negativity, Tonelli's theorem yields
	\[
		\Exf{\int_0^1 h(X_t^R)\,\dd t} = \int_0^1 \Exf{h(X_t^R)}\,\dd t\,.
	\]
	Consequently
	\[
		\Exf{\sup_{t\in[0,1]}X_t} \le \inf_{H\in\mathcal H_0} H^{-1}\left( r\int_0^1 \Exf{h(X_t^R)}\,\dd t \right)\,.
	\]

	It remains to bound $\Exf{h(X_t^R)}$ uniformly in $t$. Because $t\mapsto \Exf{f_Z(t)}$ is continuous,
	we have
	\[
		\Exf{f^R_Z(t)} = \Exf{f_Z(t)}=t\,,
	\]
	for all $t$.  Lemma~\ref{lem:RL} gives
	\[
		\eta_{f^R}\ge \eta_f=r^2\,.
	\]
	Applying McDiarmid's inequality to the fixed-parameter random variable
	$f_Z^R(t)$ gives, for every $\varepsilon>0$ and each fixed $t$,
	\[
		\mathbb P\{X_t^R\ge\varepsilon\} \le 2\e^{-2\varepsilon^2}\,.
	\]

	For fixed $t$, the layer-cake representation with respect to the
	Lebesgue--Stieltjes measure $\dd h$ gives
	\[
		\Exf{h(X_t^R)} = h(0) + \int_{(0,\infty)} \mathbb P\{X_t^R>\varepsilon\}\,\dd h(\varepsilon)\,.
	\]
	Using the preceding tail bound,
	\[
		\Exf{h(X_t^R)} \le h(0) + \int_{(0,\infty)} 2\e^{-2\varepsilon^2}\,\dd h(\varepsilon)\,.
	\]
	For $x<\infty$, integration by parts gives
	\[
		\int_{(0,x]} 2\e^{-2\varepsilon^2}\,\dd h(\varepsilon) = 2h(x)\e^{-2x^2} -2h(0) +
		\int_0^x 8\varepsilon h(\varepsilon)\e^{-2\varepsilon^2} \,\dd\varepsilon \,.
	\]
	Letting $x\to\infty$ and using $H\in\mathcal H_0$, namely
	\[
		h(x)\e^{-2x^2}\to 0,
	\]
	we obtain
	\[
		\int_{(0,\infty)} 2\e^{-2\varepsilon^2}\,\dd h(\varepsilon) = -2h(0) + \int_0^\infty 8\varepsilon h(\varepsilon)\e^{-2\varepsilon^2}
		\,\dd\varepsilon\,,
	\]
	with the identity understood in the extended sense if the final integral is
	infinite. Therefore
	\[
		\Exf{h(X_t^R)} \le -h(0) + \int_0^\infty 8\varepsilon h(\varepsilon)\e^{-2\varepsilon^2} \,\dd\varepsilon
		\le \int_0^\infty 8\varepsilon h(\varepsilon)\e^{-2\varepsilon^2} \,\dd\varepsilon,
	\]
	because $h(0)\ge0$. This bound is uniform in $t$.

	Substituting back into the Jensen--Tonelli estimate and using
	$\int_0^1\dd t=1$ yields
	\[
		\Exf{\sup_{t\in[0,1]}X_t} \le
		\inf_{H\in\mathcal H_0} H^{-1}\left( r \int_0^\infty 8\varepsilon h(\varepsilon)\e^{-2\varepsilon^2} \,\dd\varepsilon \right)\,.
	\]
	This is the claimed inequality.
\end{proof}

\subsection{The upper envelope}

For each $x>0$, let
\[
	\beta(x) := \inf_{H\in\mathcal{H}_0}
	H^{-1}\left( \sqrt{x} \int_0^\infty 8\varepsilon\, h(\varepsilon)
	\e^{-2\varepsilon^2} \,\dd\varepsilon \right)\,.
\]
By Lemma~\ref{lem:random-supremum}, any $f_Z \in \mathcal{M}_{\mathsf{id}}([0,1])$
with finite $\eta_f>0$ satisfies
\[
	0 \le \sqrt{\eta_f}\, \Exf{ \sup_{t\in[0,1]} |f_Z(t)-t| } \le
	\beta(\eta_f)\,.
\]
Define the envelope
\[
	\mathcal{E}(x) := \sqrt{\frac{2}{\ln 2}\beta(x)^2 +1} \,,
\]
for $x>0$. We obtain the following bound.

\begin{lemma}\label{prop:upper-envelope}
	For every $x>0$,
	\[
		\sqrt{U_{\mathcal{M}_{\mathsf{id}}}(x)} \le \mathcal{E}(x)\,.
	\]
\end{lemma}

\begin{proof}
	Fix $f_Z\in\mathcal{M}_{\mathsf{id}}([0,1])$ and set $x=\eta_f$. We establish that
	for all $\varepsilon>0$,
	\[
		\mathbb{P}\left\{
		\sqrt{\eta_f}\sup_{t\in[0,1]}|f_Z(t)-t| \ge
		\varepsilon\mathcal{E}(\eta_f)
		\right\}
		\le 2\e^{-2\varepsilon^2}.
	\]
	If $\varepsilon\le \sqrt{(\ln 2)/2}$, then $2\e^{-2\varepsilon^2}\ge 1$,
	and the inequality holds trivially.
	Assume $\varepsilon>\sqrt{(\ln 2)/2}$ and set $a := \sqrt{(\ln 2)/2}$. The
	$\ell_2$-norm identities are
	\[
		\mathcal{E}(\eta_f) =
		\left\| \left(1, \frac{\beta(\eta_f)}{a}\right) \right\|_2,
		\qquad
		\varepsilon =
		\left\| \left(\sqrt{\varepsilon^2-a^2}, a\right) \right\|_2 \,.
	\]
	By the Cauchy--Schwarz inequality,
	\begin{align*}
		\varepsilon\mathcal{E}(\eta_f)
		 & \ge \left\langle \left(\sqrt{\varepsilon^2-a^2}, a\right),
		\left(1, \frac{\beta(\eta_f)}{a}\right) \right\rangle         \\
		 & = \sqrt{\varepsilon^2-a^2} + \beta(\eta_f) \,.
	\end{align*}
	Define the maximal deviation
	\[
		X(Z) := \sup_{t\in[0,1]}|f_Z(t)-t| \,.
	\]
	Applying Corollary~\ref{cor:lattice_mcdiarmid} with the deterministic shift
	$\ds(t) = t$ provides the centred upper tail concentration
	\[
		\mathbb{P}\left\{
		\sqrt{\eta_f}\left( X(Z) - \Exf{X(Z)} \right) \ge v
		\right\}
		\le \e^{-2v^2}, \qquad v>0\,.
	\]
	Substituting the expectation bound
	$\sqrt{\eta_f}\Exf{X(Z)} \le \beta(\eta_f)$ of
	Lemma~\ref{lem:random-supremum}, we have
	\[
		\mathbb{P}\left\{ \sqrt{\eta_f}X(Z) \ge \beta(\eta_f)+v \right\}
		\le \e^{-2v^2}, \qquad v>0.
	\]
	Evaluating at $v = \sqrt{\varepsilon^2-a^2}$,
	\begin{align*}
		\mathbb{P}\left\{
		\sqrt{\eta_f}X(Z) \ge \varepsilon\mathcal{E}(\eta_f)
		\right\}
		 & \le \mathbb{P}\left\{
		\sqrt{\eta_f}X(Z) \ge \beta(\eta_f)+\sqrt{\varepsilon^2-a^2}
		\right\}                          \\
		 & \le \e^{-2(\varepsilon^2-a^2)}
		= \e^{-2\varepsilon^2} \e^{2a^2}
		= 2\e^{-2\varepsilon^2},
	\end{align*}
	where $\e^{2a^2}=\e^{\ln2}=2$. Hence $\mathcal{E}(x)^2$ is
	admissible, in the sense defined in Section~\ref{sec:equivalences}, for every
	$f_Z\in\mathcal{M}_{\mathsf{id}}([0,1])$ with $\eta_f=x$, so
	$\mathcal{E}(x)^2\ge u_{f_Z}$; taking the supremum, we have
	\mbox{$\sqrt{U_{\mathcal{M}_{\mathsf{id}}}(x)}\le\mathcal{E}(x)$}.
\end{proof}

\subsection{Explicit evaluation via exponential gauges}

We extract an explicit analytical envelope from $\beta(x)$. For $p>0$,
define the Young function $H_p(y) = \e^{py} - 1$, with gauge
$h_p(\varepsilon) = p \e^{p\varepsilon}$ and inverse
$H_p^{-1}(z) = \frac{1}{p}\ln(1+z)$. Because $p\e^{p\varepsilon} =
	o(\e^{2\varepsilon^2})$, we have $H_p \in \mathcal{H}_0$.
We establish the following envelope result.

\begin{proposition}\label{prop:exp-gauge}
	For any $p>0$ and $x>0$, the absolute expected supremum envelope satisfies:
	\[
		\beta(x) \le \psi(x,p) :=
		\frac{1}{p} \ln\left( 1 + 2p\sqrt{x} + p^2 \sqrt{2\pi x} \, \e^{p^2/8}
		\right) \,.
	\]
	In particular, for all $x, p>0$,
	\[
		\sqrt{U_{\mathcal{M}}(x)} \le \mathcal{E}(x,p) :=
		\sqrt{\frac{2}{\ln 2} \psi(x,p)^2 + 1}\,.
	\]
\end{proposition}

\begin{proof}
	We upper-bound the infimum in the definition of $\beta(x)$ by evaluating the
	exact integral for the specific gauge $h_p(\varepsilon) = p
		\e^{p\varepsilon}$:
	\[
		I_p := \int_0^\infty 8\varepsilon\, h_p(\varepsilon)
		\e^{-2\varepsilon^2} \,\dd \varepsilon
		= \int_0^\infty 8\varepsilon p \e^{p\varepsilon - 2\varepsilon^2}
		\,\dd\varepsilon \,.
	\]
	To evaluate this, we apply integration by parts. Let $u = 2p
		\e^{p\varepsilon}$ and $\mathrm{d} v = 4\varepsilon
		\e^{-2\varepsilon^2}\,\mathrm{d}\varepsilon$, which integrates
	to $v = -\e^{-2\varepsilon^2}$. This implies:
	\begin{align*}
		I_p
		 & = \left[ -2p \e^{p\varepsilon - 2\varepsilon^2} \right]_0^\infty
		- \int_0^\infty \left(2p^2 \e^{p\varepsilon}\right)
		\left(-\e^{-2\varepsilon^2}\right) \,\mathrm{d}\varepsilon          \\
		 & = 2p + 2p^2 \int_0^\infty \e^{p\varepsilon - 2\varepsilon^2}
		\,\mathrm{d}\varepsilon \,.
	\end{align*}
	We isolate the Gaussian core by completing the square in the exponent:
	\[
		-2\varepsilon^2 + p\varepsilon = -2\left(\varepsilon - \frac{p}{4}\right)^2
		+ \frac{p^2}{8}\,.
	\]
	Substituting this back into the integral and factoring out
	the exponential, we have
	\[
		I_p = 2p + 2p^2 \e^{p^2/8} \int_0^\infty
		\e^{-2(\varepsilon - p/4)^2} \,\mathrm{d}\varepsilon \,.
	\]
	We upper-bound the remaining shifted Gaussian integral by extending
	its lower limit to $-\infty$, which evaluates to $\sqrt{\pi/2}$.
	Substituting this, we have
	\[
		I_p \le 2p + 2p^2 \sqrt{\frac{\pi}{2}} \e^{p^2/8}
		= 2p + p^2\sqrt{2\pi} \, \e^{p^2/8} \,.
	\]
	Because $\beta(x)$ is defined as an infimum, it is bounded by the evaluation
	at $H_p$, so that
	\[
		\beta(x) \le H_p^{-1}(\sqrt{x} I_p) =
		\frac{1}{p} \ln(1 + \sqrt{x} I_p)\,.
	\]

	The ``in particular'' claim follows because
	$\sqrt{U_{\mathcal{M}_{\mathsf{id}}}(x)}\le\mathcal{E}(x)\le\mathcal{E}(x,p)$
	by Lemma~\ref{prop:upper-envelope}, and the uniformisation
	Theorem~\ref{thm:unif} transports the bound to $U_{\mathcal{M}}$, the
	envelope $\mathcal{E}(\cdot,p)$ being continuous.
\end{proof}

\subsubsection{Large-value tuning}

Tuning $p(x)=2\sqrt{\ln x}$ in Proposition~\ref{prop:exp-gauge},
and writing $y=\ln x$, the choice $p=2\sqrt{y}$ has $p^2/8=y/2$, so that
$\psi(\e^y,2\sqrt{y})=\tfrac{\sqrt{y}}{2}\bigl(1+\w(y)\bigr)$ and hence
\begin{equation}\label{eq:lemasymp}
	\sqrt{U_{\mathcal{M}}(\e^y)}\le \mathcal{E}(\e^y,2\sqrt{y})
	= \sqrt{\frac{y}{2\ln 2}}\,
	\sqrt{(1+\w(y))^2 + \frac{2\ln 2}{y}}\,,
\end{equation}
where
\[
	\w(y) := \frac1y \ln\left(
	4\sqrt{2\pi}\,y + 4\sqrt{y}\,\e^{-y/2}
	+ \e^{-y} \right) \,.
\]

The proof of the following computational lemma is deferred to the Appendix.

\begin{lemma}\label{lem:scalar-envelope}
	With $\w$ as in \eqref{eq:lemasymp}, the bound
	\[
		\sqrt{\frac{y}{2\ln 2}}
		\left[ \sqrt{ (1+\w(y))^2 + \frac{2\ln 2}{y} } - 1 \right]
		\le \frac{2\ln(\e+y)}{\sqrt{y}}
	\]
	holds for all $y \ge y_0$, where $y_0$ is the unique solution of
	$\ln(\e+y)/\sqrt{y} = 1$, $y \ge 1$.
\end{lemma}

We obtain the needed bound beyond the threshold, for $y\ge y_0$.

\begin{corollary}\label{cor:asymp-envelope}
	For every $x$ with $\ln x\ge y_0$,
	\[
		\sqrt{U_{\mathcal{M}}(x)}-\sqrt{\log_4 x}
		\le \frac{2\ln(\e+\ln x)}{\sqrt{\ln x}}\,.
	\]
\end{corollary}

\begin{proof}
	Write $y=\ln x$, so $\sqrt{\log_4 x}=\sqrt{y/(2\ln2)}$. Subtracting this
	leading term from \eqref{eq:lemasymp}, we have
	\[
		\sqrt{U_{\mathcal{M}}(\e^y)}-\sqrt{\tfrac{y}{2\ln2}}
		\le \sqrt{\tfrac{y}{2\ln2}}
		\left[\sqrt{(1+\w(y))^2+\tfrac{2\ln2}{y}}-1\right] \,,
	\]
	which by Lemma~\ref{lem:scalar-envelope} is at most
	$2\ln(\e+y)/\sqrt{y}$ for $y\ge y_0$.
\end{proof}

\subsubsection{Small-value tuning}

For small $x$, let $p=5$ in Proposition~\ref{prop:exp-gauge}, chosen heuristically.
The substitution provides $\psi(x,5)=\tfrac15\ln(1+K\sqrt{x})$ with
$K=10+25\sqrt{2\pi}\,\e^{25/8}$, hence
\begin{equation}\label{eq:lemsmall}
	\sqrt{U_{\mathcal{M}}(x)}\le \mathcal{E}(x,5)
	= \sqrt{\frac{2}{\ln 2}\,\psi(x,5)^2 + 1} \,.
\end{equation}

The following is a direct calculation over the range $[2.5, \e^{y_0})$.

\begin{lemma}\label{lem:small-envelope}
	For every $x$ with $2.5\le x$ and $\ln x<y_0$,
	\[
		\sqrt{U_{\mathcal{M}}(x)}-\sqrt{\log_4 x} < 2 \,.
	\]
\end{lemma}

\begin{proof}
	By \eqref{eq:lemsmall},
	\[
		\sqrt{U_{\mathcal{M}}(x)}-\sqrt{\log_4 x}
		\le \mathcal{E}(x,5)-\sqrt{\log_4 x}=:S(x) \,.
	\]
	Because $\mathcal{E}(x,5)\ge\sqrt{\tfrac{2}{\ln2}}\,\psi(x,5)$, on
	$2.5\le x<\e^{y_0}$
	\[
		\frac{\mathrm d}{\mathrm dx}\,\mathcal{E}(x,5)
		= \frac{(2/\ln2)\,\psi(x,5)\,\psi'(x,5)}{\mathcal{E}(x,5)}
		\le \sqrt{\tfrac{2}{\ln2}}\,\psi'(x,5)
		< \frac{\mathrm d}{\mathrm dx}\sqrt{\log_4 x} \,,
	\]
	so $S$ is decreasing and attains its maximum at $x=2.5$:
	\[
		S(2.5)=\mathcal{E}(2.5,5)-\sqrt{\log_4 2.5} < 2,  \qquad(1.996547\ldots)\,.
	\]
	This establishes the result.
\end{proof}

\subsubsection{Microrange}

We also know the trivial bound of Lemma~\ref{lem:u-upper}, which holds for all
classes of $[0,1]$-valued randomised functions:
\[
	\sqrt{U_\mathcal{M}(x)} \le \sqrt{\frac{2}{\ln2}\,x},
	\qquad \text{for all } x\ge 0 \,.
\]
It is the sharpest of the three on the microrange $[1,2.5)$.

\begin{lemma}\label{lem:tiny-envelope}
	For every $x$ with $1\le x < 2.5$,
	\[
		\sqrt{U_{\mathcal{M}}(x)}-\sqrt{\log_4 x} < 2 \,.
	\]
\end{lemma}

\begin{proof}
	By the trivial bound of Lemma~\ref{lem:u-upper},
	\[
		\sqrt{U_{\mathcal{M}}(x)}-\sqrt{\log_4 x}\le\gs(x):=
		\sqrt{\tfrac{2}{\ln2}\,x}-\sqrt{\log_4 x}\,.
	\]
	A direct verification gives $\gs(x)\le \gs(2.5)\le 1.8729<2$ on $[1,2.5]$.
\end{proof}

\begin{proof}[Proof of Theorem~\ref{thm:upper_bound_main}]
	Put $y=\ln x$, so that $\sqrt{\log_4 x}=\sqrt{y/(2\ln2)}$, and let
	$x_0=\e^{y_0}$. It suffices to show
	\[
		\sqrt{U_{\mathcal{M}}(x)}-\sqrt{\log_4 x}
		\le 2\min\!\left\{1,\frac{\ln(\e+y)}{\sqrt{y}}\right\}
	\]
	for all $x>1$.
	The function $g(y):=\ln(\e+y)/\sqrt{y}$ is strictly decreasing on
	$(0,\infty)$, because $g'(y)<0$ is equivalent to the elementary inequality
	$2y/(\e+y)<\ln(\e+y)$; and $g(y_0)=1$, so $g\le1$ for $y\ge y_0$
	and $g\ge1$ for $0<y\le y_0$. Hence
	\[
		2\min\!\left\{1,\frac{\ln(\e+y)}{\sqrt{y}}\right\}
		=\begin{cases}
			\dfrac{2\ln(\e+y)}{\sqrt{y}}, & y\ge y_0 \quad(x\ge x_0),   \\[8pt]
			2,                            & 0<y< y_0 \quad(1<x<x_0) \,.
		\end{cases}
	\]
	\emph{Case $x\ge x_0$.} By Corollary~\ref{cor:asymp-envelope},
	\[
		\sqrt{U_{\mathcal{M}}(x)}-\sqrt{\log_4 x}
		\le \frac{2\ln(\e+\ln x)}{\sqrt{\ln x}}
		= 2\min\!\left\{1,\frac{\ln(\e+y)}{\sqrt{y}}\right\} \,.
	\]
	\emph{Case $1<x<x_0$.} Here the right-hand side is $2$. The ranges $[1,2.5)$
	and $[2.5,x_0)$ cover $(1,x_0)$: on $[1,2.5)$,
	Lemma~\ref{lem:tiny-envelope} provides
	$\sqrt{U_{\mathcal{M}}(x)}-\sqrt{\log_4 x}<2$, and on $[2.5,x_0)$,
	Lemma~\ref{lem:small-envelope} provides the same bound. Thus
	$\sqrt{U_{\mathcal{M}}(x)}-\sqrt{\log_4 x}<2$ throughout $(1,x_0)$.
	Combining the two cases proves the bound for every $x>1$.
\end{proof}

\section{Proof of Theorem~\ref{thm:main}}\label{sec:proof_main}

Combining Theorem~\ref{thm:lower_bound_main} and Theorem~\ref{thm:upper_bound_main}
gives, for every $x > 1$, the two-sided bound
\begin{equation}\label{eq:asy}
	\begin{aligned}
		-2\min\!\left\{1, \, \frac{2\ln(\e+\ln x)}{\sqrt{\ln x}}\right\}
		 & \le \sqrt{U_{\mathcal{M}}(x)} - \sqrt{\log_4 x}                        \\
		 & \le 2\min\!\left\{1, \, \frac{\ln(\e+\ln x)}{\sqrt{\ln x}}\right\} \,,
	\end{aligned}
\end{equation}
which gives the symmetric inequality of Theorem~\ref{thm:main}.

\section{Applications}\label{sec:applications}

We open with an empirical distribution function for stratified clustered
data that exhibits the lower bound (Theorem~\ref{thm:lower_bound_main}) as sharp.
We then give an application of the upper bound
to non-monotone bounded processes: an analogue of the smooth Gaussian
setting of Diebolt and Posse~\cite{diebolt1996density} with the
$\ell_1$ arc length in place of the $\ell_2$ arc length used there.
Although this is one instance of the randomised signed measures of
Proposition~\ref{prop:maharam_cumulative}, we treat it here natively.

\subsection{Empirical distribution function with stratified data}\label{sec:EDFcluster}

We present an empirical distribution function that asymptotically
saturates the main theorem. It arose in the analysis of a randomised
field experiment with clustered structure, the cluster sizes growing
faster than the number of clusters, and motivated the present study.
The classical uniform law of large numbers does not in general apply
to such data (cf.~\cite{canay2018wild}); the stratification below
restores it.

Let
\[
	X_{i,j}, \quad 1 \le i \le n, ~1 \le j \le m_n \,,
\]
be real-valued random variables, the blocks $(X_{i,1}, \dots, X_{i,m_n})$
being independent across $i$, and write
\begin{equation}\label{eq:F}
	F_{nm}(t) = \frac{1}{n m_n} \sum_{i=1}^n \sum_{j=1}^{m_n}
	\mathbf{1}_{\{X_{i,j} \le t\}} \,.
\end{equation}
Suppose the supports are disjoint and ordered within each block:
\begin{equation}\label{eq:stratification}
	\operatorname*{ess\,sup} X_{i,j} \le \operatorname*{ess\,inf} X_{i,j+1} \,.
\end{equation}
Any threshold $t$ lies in at most one of the intervals
$[\operatorname*{ess\,inf} X_{i,j}, \operatorname*{ess\,sup} X_{i,j}]$,
and resampling the $i$-th block alters the inner sum
$\sum_{j} \mathbf{1}_{\{X_{i,j} \le t\}}$ by at most $1$. Hence
$\eta_{F_{nm}} \ge n m_n^2$, with equality holding in the non-degenerate case.

An instance of $F_{nm}$ realises the lower-bound construction of
Section~\ref{sec:lower_bound}. Take $\mathcal{Z} = \{-1, +1\}^n$ and a
non-empty $\Theta \subseteq \mathcal{Z}$ with $|\Theta| = m_n$. Arrange
$\Theta = [\theta_{ij}]_{i,j=1}^{n, m_n}$ as a matrix, with columns
$\theta_j$. Draw $Z$ uniformly on $\mathcal{Z}$ and set
\[
	X_{i,j} = Z_i \theta_{ij} + 4j \,.
\]
The supports $\{4j-1, 4j+1\}$ are disjoint and ordered, so
\eqref{eq:stratification} holds. For $t \in [4j-1, 4j+1)$,
\[
	F_{nm}(t) - \Exf{F_{nm}(t)}
	= -\frac{1}{2 n m_n} \langle Z, \theta_j \rangle \,,
\]
and the centred difference vanishes elsewhere. The right-hand side is
the inner product $f$ of~\eqref{eq:f-extreme}, scaled by
$-(2 n m_n)^{-1}$; by Lemma~\ref{lem:elementary}, $F_{nm}$ is equivalent
to $f_Z$ on $\Theta$, with \mbox{$\eta_{F_{nm}} = n m_n^2$} as before.

Consider the regime
\begin{equation*}
	\frac{\ln n}{\ln m_n} \to 0 \,,
\end{equation*}
in which $\ln n / \ln \eta_{F_{nm}} \to 0$ and
$\ln m_n / \ln \eta_{F_{nm}} \to 1/2$. The tuned constellation construction of
Section~\ref{subsec:calibrated}, applied with
$\ln m_n \sim \tfrac{5}{2}\sqrt{n}$, gives
\[
	\lim_{n \to \infty} \frac{u_{F_{nm}}}{\ln \eta_{F_{nm}}}
	= \frac{1}{2 \ln 2} \,,
\]
the constant of the main theorem; with this growth we realise
finite-sample lower bounds in the proof of
Theorem~\ref{thm:lower_bound_main}.

This cluster-growth regime, attaining the optimal normalisation
$\sqrt{\eta_{F_{nm}}/\log_4\eta_{F_{nm}}}$ of the main theorem,
is natural in cluster sampling: under the Ewens sampling formula, for example,
the $k$-th largest cluster grows linearly while the number of clusters grows
logarithmically.

The intermediate regime,
\[
	\frac{\ln m_n}{\ln n} \to c \in (0,\infty) \,,
\]
for which the bracketing normalisation $\sqrt{n}$ is informative,
is incompatible with cluster sampling models admitting a Kingman paintbox
representation~\cite{kingman1978representation,pitman1995exchangeable}
and is realised by non-exchangeable microclustering
sampling~\cite{dibenedetto2021nonexchangeable}.

The reverse regime, in which $n$ outpaces $m_n$,
\[
	\frac{\ln m_n}{\ln n} \to 0 \,,
\]
for which the bracketing normalisation $\sqrt{n}$ is tight,
is unattainable with microclustering sampling, and for exchangeable
sampling is realised only in the trivial case
\[
	m_n = 1 \,,
\]
recovering Massart's framework of scalar i.i.d.\ empirical
distribution functions.

\begin{remark}\label{rem:bracketing_loose}
	The cluster empirical distribution $F_{nm}$ is an average of $n$ independent
	monotone summands. For the identically distributed case, the bracketing
	estimates~\cite{vanderVaartWellner1996} normalise by $\sqrt{n}$, where $n$ is
	the number of independent summands, ensuring
	\[
		\Exf{\sup_t \left|F_{nm}(t) - \Exf{F_{nm}(t)}\right|} \le C/\sqrt{n}\,,
	\]
	for a
	universal constant $C$. The optimal normalisation of Theorem~\ref{thm:main}
	instead uses the effective sample size $\eta_{F_{nm}} = n m_n^2$, normalising
	by $\sqrt{\eta_{F_{nm}}/U_{\mathcal{M}}(\eta_{F_{nm}})} \sim
		\sqrt{\eta_{F_{nm}}/\log_4 \eta_{F_{nm}}}$. The bracketing normalisation is,
	therefore, loose by the factor
	\[
		\frac{\sqrt{\eta_{F_{nm}}/\log_4 \eta_{F_{nm}}}}{\sqrt{n}}
		= \frac{m_n}{\sqrt{\log_4 \eta_{F_{nm}}}}
		\sim \frac{m_n}{\sqrt{2 \log_4 m_n}}\,,
	\]
	in the regime $\ln n / \ln m_n \to 0$. In a regime realising the lower-bound construction
	(Section~\ref{sec:lower_bound}), $\ln m_n \sim \tfrac{5}{2}\sqrt{n}$, so
	bracketing's $\sqrt{n}$ normalisation is exponentially loose at the scale of the
	effective sample size $\sqrt{\eta_{F_{nm}}} = \sqrt{n}\, m_n$.
\end{remark}

\subsection{Arc length of smooth bounded processes}\label{subsec:arc_length}

We turn to an application of the upper bound
(Theorem~\ref{thm:upper_bound_main}), and present the class of smooth
bounded processes whose suprema are controlled by the $\ell_1$ arc length of
the coordinate functions. This is the bounded distribution-free counterpart of a
classical Gaussian extreme-value setting.

Consider a randomised function
\[
	X_Z(t) = \sum_{j=1}^m g_j(Z) \psi_j(t),
	\qquad t \in [0, T] \,,
\]
where $g_j \colon \mathcal{Z} \to [-1, 1]$ are measurable, and
$\psi_j \colon [0, T] \to \mathbb{R}$ are continuously
differentiable with $\sum_{j=1}^m |\psi_j(t)| = 1$ for all
$t \in [0, T]$.

Assume $\Exf{X_Z(t)} = 0$ for all $t \in [0, T]$.
Define the $\ell_1$ arc length
\[
	\mathrm{Arc}_{\ell_1}(t) := \int_0^t \sum_{j=1}^m |\psi'_j(s)| \, \mathrm{d}s,
	\qquad \mathrm{Arc}_{\ell_1} := \mathrm{Arc}_{\ell_1}(T) \,.
\]

\begin{corollary}\label{cor:differentiable_paths}
	If $1< \eta_X<\infty$, then for all $\varepsilon > 0$,
	\[
		\mathbb{P}\!\left\{
		\frac{\sqrt{\eta_X}}{\mathcal{R}\bigl((2 + \mathrm{Arc}_{\ell_1})^2 \eta_X\bigr)}
		\sup_{t \in [0, T]} |X_Z(t)| \ge \varepsilon
		\right\}
		\le 2 \e^{-2\varepsilon^2} \,,
	\]
	where
	\[
		\mathcal{R}(x) := \sqrt{\log_4 x} +
		2\min\!\left\{1, \frac{\ln(\e+\ln x)}{\sqrt{\ln x}}\right\},
		\qquad  \text{for all $x> 1$} \,.
	\]
\end{corollary}

\begin{proof}
	The bound $|X_Z(t)| \le 1$ holds by the normalisation
	$\sum_{j=1}^m |\psi_j(t)| = 1$, and the local-derivative bound
	$\kappa(t) := \sum_{j=1}^m |\psi'_j(t)|$ satisfies
	$|X'_Z(t)| \le \kappa(t)$ for every realisation. The drifted process
	\[
		Y_Z(t) := X_Z(t) + \mathrm{Arc}_{\ell_1}(t) + 1
	\]
	has derivative $Y'_Z(t) = X'_Z(t) + \kappa(t) \ge 0$, so for every
	realisation $Z = z$ the path $Y_z$ is non-decreasing in $t$, and
	$Y_Z(t) \in [0,\, 2 + \mathrm{Arc}_{\ell_1}]$. Normalise to define
	\[
		f_Z(t) := \frac{Y_Z(t)}{2 + \mathrm{Arc}_{\ell_1}} \,,
	\]
	giving $f_Z \in \mathcal{M}$. Applying Lemma~\ref{lem:elementary} twice,
	first with $a = 1$ and deterministic shift
	\[
		\ds(t) = \mathrm{Arc}_{\ell_1}(t) + 1\,,
	\]
	second with
	$a = (2 + \mathrm{Arc}_{\ell_1})^{-1}$ and $\ds \equiv 0$, shows that $f_Z$
	is equivalent to $X_Z$ with $\eta_{f} = (2 + \mathrm{Arc}_{\ell_1})^2
		\eta_X>1$, so Theorem~\ref{thm:upper_bound_main} applied to $f_Z$ and
	Lemma~\ref{cor:U} give the result.
\end{proof}

\appendix
\section{Preliminaries}

\subsection{Characterisation of the class by randomised signed-measure
	processes}

\begin{proof}[Proof of Proposition~\ref{prop:maharam_cumulative}]
	Let $F_Z$ be the given regular cumulative function. Use the separating
	measure $\lambda$ to scalarise the chain. Set
	\[
		T:=\{\lambda(E_\theta):\theta\in\Theta\}\subseteq\mathbb R\,.
	\]
	Because the chain is totally ordered by inclusion and separated by
	$\lambda$-positive symmetric differences, the map
	\[
		\theta\longmapsto \lambda(E_\theta)
	\]
	is injective. Hence, for each $t\in T$, there is a unique
	$\theta(t)\in\Theta$ such that
	\[
		\lambda(E_{\theta(t)})=t\,.
	\]
	Write
	\[
		E_t:=E_{\theta(t)}, \qquad t\in T\,.
	\]
	For $s,t\in T$,
	\[
		s\ge t \quad\Longrightarrow\quad E_t\subseteq E_s\,.
	\]

	Define
	\[
		g\colon\mathcal Z\times T\to\mathbb R, \qquad
		g(z,t):=\mu_z(E_t)=\langle \mu_z,\mathbf 1_{E_t}\rangle\,.
	\]
	The randomised function $g_Z\colon T\to\mathbb R$ is the scalar-indexed
	reparametrisation of $F_Z$. Hence $g_Z$ is equivalent to $F_Z$, with
	$\eta_g=\eta_F$.

	Set
	\[
		\nu:=\bigwedge_{z\in\mathcal Z}\mu_z\,,
	\]
	which exists by Dedekind completeness, because the family is order
	bounded from below. We have $\mu_z-\nu\ge0$ as a measure for every
	$z\in\mathcal Z$. Define
	\[
		h_Z(t)
		:= \frac{\langle \mu_Z-\nu,\mathbf 1_{E_t}\rangle}{1+\|\nu\|}
		= \frac{g_Z(t)-\nu(E_t)}{1+\|\nu\|}, \qquad t\in T\,.
	\]
	Because $\mu_z-\nu$ is positive and $s\ge t$ implies
	$E_t\subseteq E_s$, the function $h_z$ is non-decreasing in the scalar
	parameter $t$ for each $z$. Moreover,
	\[
		0\le (\mu_z-\nu)(E_t)
		\le \|\mu_z-\nu\|
		\le \|\mu_z\|+\|\nu\|
		\le 1+\|\nu\|\,,
	\]
	so $h_Z(t)\in[0,1]$. Thus, $h_Z$ is a $[0,1]$-valued monotone randomised
	function on $T$.

	The randomised function $h_Z$ is equivalent to $g_Z$, and hence to
	$F_Z$, with
	\[
		\eta_h
		= \bigl(1+\|\nu\|\bigr)^2\,\eta_g
		= \bigl(1+\|\nu\|\bigr)^2\,\eta_F\,.
	\]
	To obtain an element of $\mathcal M$ indexed by all of $\mathbb R$,
	extend $h_Z$ by
	\[
		f(z,r)
		:= \sup\{h(z,t):t\in T,\ t\le r\},
		\qquad (z,r)\in\mathcal Z\times\mathbb R\,,
	\]
	with the convention that the supremum over the empty set is $0$.
	The sections $z\mapsto f_z(t)$ are measurable due to monotonicity.
	So we have $f_Z\in\mathcal M$. Further, this monotone extension $f_Z$ is
	equivalent to $h_Z$ and leaves the effective sample size unchanged,
	\[
		\eta_f = \eta_h = \bigl(1+\|\nu\|\bigr)^2\,\eta_F\,,
	\]
	and $f_Z$ is equivalent to $F_Z$.

	Conversely, let $\preceq$ be the lexicographic order on $\mathbb R^3$
	and set
	\[
		E_t^\ast := \{\,x\in\mathbb R^3:\ x\preceq (0,t,1)\,\},
		\qquad t\in\mathbb R\,.
	\]
	Equivalently,
	\[
		E_t^\ast
		= (-\infty,0)\times\mathbb R^2
		\;\cup\; \{0\}\times(-\infty,t)\times\mathbb R
		\;\cup\; \{0\}\times\{t\}\times(-\infty,1]\,.
	\]
	The family $(E_t^\ast)_{t\in\mathbb R}$ is an increasing Borel chain.
	Moreover, it is a regular chain. Indeed, let $\varphi$ be any strictly
	positive integrable density on $\mathbb R$, and define a finite positive
	Borel measure $\lambda$ on $\mathbb R^3$ by
	\[
		\lambda(A)
		:= \int_{\mathbb R} \mathbf 1_A(0,u,1)\,\varphi(u)\,\dd u,
		\qquad A\in\mathcal B(\mathbb R^3)\,.
	\]
	Equivalently, $\lambda$ is the push-forward of $\varphi(u)\,\dd u$ under
	the map $u\mapsto(0,u,1)$. Thus, $\lambda$ is a finite positive Borel
	measure on $\mathbb R^3$, and for $s<t$,
	\[
		\lambda(E_t^\ast\triangle E_s^\ast) = \int_s^t \varphi(u)\,\dd u > 0\,.
	\]
	Hence,
	\[
		s\ne t
		\quad\Longrightarrow\quad
		\lambda(E_s^\ast\triangle E_t^\ast)>0\,,
	\]
	so $(E_t^\ast)_{t\in\mathbb R}$ is a regular chain.

	Fix $f_Z\in\mathcal M$. Choose a strictly increasing bijection
	$\psi\colon(0,1)\to\mathbb{R}$ and define
	\[
		\tilde f(z,t) :=
		\begin{cases}
			0,            & t\le 0\,, \\
			f(z,\psi(t)), & 0<t<1\,,  \\
			1,            & t\ge 1\,.
		\end{cases}
	\]
	By Lemma~\ref{lem:elementary}, $\tilde f_Z$ is equivalent to $f_Z$ with
	$\eta_{\tilde f}=\eta_f$. Each section $\tilde f_z$ is non-decreasing
	with values in $[0,1]$, so $\tilde f_Z\in\mathcal M$.

	For each $z\in\mathcal{Z}$ there is a positive finite Borel measure
	$\mu_z$ on $\mathbb R^3$ with $\|\mu_z\|\le 1$ such that
	\[
		\tilde f(z,t)=\mu_z(E_t^\ast), \qquad t\in\mathbb R\,.
	\]
	Indeed, push the continuous part of $\tilde f_z$ onto
	$\{0\}\times\mathbb R\times\{1\}$ via its Lebesgue--Stieltjes measure;
	place the left jump $\tilde f_z(t)-\tilde f_z^{\mathrm L}(t)$ as an atom
	at $(0,t,1)$ and the right jump $\tilde f_z^{\mathrm R}(t)-\tilde f_z(t)$
	as an atom at $(0,t,2)$. Because $\tilde f_z$ vanishes on $(-\infty,0]$,
	the increments collected by $E_t^\ast$ total $\tilde f_z(t)$, so
	$\mu_z(E_t^\ast)=\tilde f_z(t)$; the total mass is
	$\lim_{t\to\infty}\tilde f_z^{\mathrm R}(t)=1$.

	The family $\{\mu_z\}_{z\in\mathcal Z}$ is order bounded from below by
	the zero measure, so $\nu_0:=\bigwedge_{z\in\mathcal Z}\mu_z$ exists by
	Dedekind completeness and satisfies $0\le\nu_0\le\mu_z$. The family
	$\{\mu_z-\nu_0\}_{z\in\mathcal Z}$ consists of positive measures with
	$\|\mu_z-\nu_0\|\le\|\mu_z\|\le1$, and $\bigwedge_{z}(\mu_z-\nu_0)=0$. Its
	cumulative function is $\tilde f_Z-\ds$, where
	$\ds(t):=\nu_0(E_t^\ast)$ is deterministic and bounded, so it is
	equivalent to $\tilde f_Z$, and thus to $f_Z$, by another application of
	Lemma~\ref{lem:elementary}. The required pointwise measurability follows
	from the measurability of the sections of $\tilde f$.
\end{proof}

\section{The lower bound}

\subsection{Central binomial bound for constellations}\label{sec:appendix_binomial}

To prove the required binomial lower bound for a constellation, we use the bounds
\begin{equation}\label{eq:binbound}
	\frac{2^{n\hspace{0.5pt}\Hb(\lambda)}}{\sqrt{8\hspace{0.5pt}n\lambda(1-\lambda)}}
	\le \binom{n}{\lambda n} \le
	\frac{2^{n\hspace{0.5pt}\Hb(\lambda)}}{\sqrt{2\pi n\lambda(1-\lambda)}},
	\qquad 0<\lambda<1,~ n\lambda\in\mathbb{N},~ n\in\mathbb{N} \,,
\end{equation}
where $\Hb \colon [0,1] \to [0,1]$,
\[
	\Hb(x) := -x\log_2 x - (1-x)\log_2(1-x) \,,
\]
is the binary entropy function \cite[Ch.~10, Lemma~7, p.~309]{macwilliams1977theory}.
This function is symmetric about its maximiser $x=\tfrac12$ and
the exact Taylor expansion about $\tfrac{1}{2}$ has no odd terms,
\begin{equation}\label{eq:H2_taylor}
	1-\Hb\left(\frac{1}{2}-\delta\right) =
	\frac{1}{\ln 2} \sum_{j=1}^{\infty} \frac{(2\delta)^{2j}}{2j(2j-1)},
	\qquad \text{for } |\delta| < 1/2 \,.
\end{equation}

We are ready to prove the required bound.

\begin{lemma}\label{lem:true_binomial_term_general}
	For a constellation $(k, t_k, n_k, r_k, m_k)$, we have
	\[
		\frac{1}{2^{n_k}}\binom{n_k}{r_k} \ge \frac{1}{\sqrt{2k}\,t_k}
		\exp\left(-\frac{k}{12(t_k^2-1)}\right) \e^{-k/2} \,.
	\]
\end{lemma}

\begin{proof}
	In light of Lemma~\ref{lem:constellation_from_parameters}, for the constellation,
	set
	\[
		\lambda_k := \frac{r_k}{n_k} =
		\frac{k\,t_k(t_k-1)}{2k\,t_k^2}= \frac{t_k-1}{2t_k}
		= \frac{1}{2} - \frac{1}{2t_k} \,.
	\]
	Note that $1/4\le\lambda_k <1/2$, because $t_k\ge 2$.
	Thus, we apply the lower bound of \eqref{eq:binbound}:
	\[
		\binom{n_k}{r_k}= \binom{n_k}{\lambda_k n_k}
		\ge \frac{1}{\sqrt{8n_k\lambda_k(1-\lambda_k)}}\,
		2^{n_k \Hb(\lambda_k)} \,.
	\]

	Dividing by $2^{n_k}$ gives
	\begin{equation}\label{eq:brackterm}
		\frac{1}{2^{n_k}}\binom{n_k}{r_k} \ge
		\frac{1}{\sqrt{8n_k\lambda_k(1-\lambda_k)}} \, \left[2^{-n_k(1-\Hb(\lambda_k))}\right] \,.
	\end{equation}

	Now
	\[
		\lambda_k(1-\lambda_k) = \left(\frac{1}{2}-\frac{1}{2t_k}\right)
		\left(\frac{1}{2}+\frac{1}{2t_k}\right) = \frac{1}{4}-\frac{1}{4t_k^2}
		\le \frac{1}{4} \,,
	\]
	and therefore
	\begin{equation}\label{eq:firstest}
		\frac{1}{\sqrt{8n_k\lambda_k(1-\lambda_k)}} \ge \frac{1}{\sqrt{2n_k}}
		= \frac{1}{\sqrt{2k}\,t_k} \,.
	\end{equation}

	It remains to bound the bracketed term in \eqref{eq:brackterm}.
	Writing $\delta:=\frac{1}{2t_k}$ and noting that
	\mbox{$0<\delta \le 1/4$}, we have
	\[
		\lambda_k = \frac{1}{2} - \delta \,.
	\]
	By the exact Taylor expansion in \eqref{eq:H2_taylor}, we have
	\[
		1-\Hb\left(\frac{1}{2}-\delta\right) =
		\frac{1}{\ln 2} \sum_{j=1}^{\infty} \frac{(2\delta)^{2j}}{2j(2j-1)} \,.
	\]
	Because $2\delta=1/t_k$, this becomes
	\[
		1-\Hb(\lambda_k) = \frac{1}{\ln 2} \sum_{j=1}^{\infty}
		\frac{t_k^{-2j}}{2j(2j-1)} \,.
	\]
	Separating the first term and bounding the remainder by a geometric series,
	\begin{align*}
		1-\Hb(\lambda_k)
		 & \le \frac{1}{\ln 2} \left[ \frac{1}{2t_k^2} +
		\frac{1}{12}\sum_{j=2}^{\infty}t_k^{-2j} \right] \\
		 & = \frac{1}{\ln 2} \left[
		\frac{1}{2t_k^2} + \frac{1}{12t_k^4}\sum_{m=0}^{\infty}t_k^{-2m} \right] \,.
	\end{align*}
	Hence,
	\[
		1-\Hb(\lambda_k) \le \frac{1}{\ln 2} \left[ \frac{1}{2t_k^2} +
			\frac{1}{12t_k^2(t_k^2-1)} \right] \,,
	\]
	and multiplying by $n_k\ln 2=k\,t_k^2\ln 2$ gives
	\[
		n_k(1-\Hb(\lambda_k))\ln 2 \le k\,t_k^2\left(\frac{1}{2t_k^2} +
		\frac{1}{12t_k^2(t_k^2-1)}\right) = \frac{k}{2} + \frac{k}{12(t_k^2-1)} \,.
	\]
	Thus
	\begin{align*}
		2^{-n_k(1-\Hb(\lambda_k))}
		 & = \exp\left(-n_k(1-\Hb(\lambda_k))\ln 2\right)          \\
		 & \ge \exp\left(-\frac{k}{2}-\frac{k}{12(t_k^2-1)}\right)
		= \exp\left(-\frac{k}{12(t_k^2-1)}\right)\e^{-k/2} \,.
	\end{align*}

	Substituting this and \eqref{eq:firstest}
	into \eqref{eq:brackterm}, we obtain
	\[
		\frac{1}{2^{n_k}}\binom{n_k}{r_k} \ge \frac{1}{\sqrt{2k}\,t_k}
		\exp\left(-\frac{k}{12(t_k^2-1)}\right) \e^{-k/2} \,.
	\]
	This proves the claim.
\end{proof}

\section{Parameter-side uniformisation}

\subsection{Right-continuous and left-continuous versions}

For a randomisation $(f,Z)$ belonging to $\mathcal M$,
let $f^R$ and $f^L$ be defined section-wise as follows:
for each $z\in\mathcal Z$ and $t\in\mathbb R$,
\[
	f^R_z(t) := \lim_{s\downarrow t} f_z(s) = \inf_{s>t} f_z(s) \,,
\]
and
\[
	f^L_z(t) := \lim_{s\uparrow t} f_z(s) = \sup_{s<t} f_z(s) \,.
\]

\begin{lemma}\label{lem:RL}
	The pairs $(f^R,Z)$ and $(f^L,Z)$ are randomisations belonging to $\mathcal M$. Moreover, for every
	$i=1,\ldots,n$,
	\[
		c_i(f^R)=c_i(f^L)\le c_i(f) \,.
	\]

	Furthermore, for every $z\in\mathcal Z$,
	\[
		\sup_{t\in\mathbb R}
		\left| f^R_z(t)-\Exf{f^R_Z(t)} \right| = \sup_{t\in\mathbb R} \left| f^L_z(t)-\Exf{f^L_Z(t)} \right|
		\le \sup_{t\in\mathbb R} \left| f_z(t)-\Exf{f_Z(t)} \right| \,.
	\]
\end{lemma}

\begin{proof}
	Because $f_z$ is non-decreasing and takes values in $[0,1]$, the limits
	defining $f^R_z(t)$ and $f^L_z(t)$ exist for every $z$ and $t$. By the density of the rationals,
	\[
		f^R_z(t) = \inf_{\substack{q\in\mathbb{Q}\\ q>t}} f_z(q),
		\qquad
		f^L_z(t) = \sup_{\substack{q\in\mathbb{Q}\\ q<t}} f_z(q) \,.
	\]
	For fixed $t$, these are respectively a countable infimum and a countable
	supremum of measurable functions of $z$. Hence, $z\mapsto f^R_z(t)$ and
	$z\mapsto f^L_z(t)$ are measurable. Moreover, for fixed $z$, the functions
	$f^R_z$ and $f^L_z$ are non-decreasing and take values in $[0,1]$. Thus,
	$(f^R,Z)$ and $(f^L,Z)$ belong to $\mathcal{M}$.

	We next record a deterministic fact. Let $a,b$ be functions in $\mathsf{M}$.
	Let $a^R,b^R,a^L,b^L$ denote their right-continuous and left-continuous versions.
	We have
	\[
		\sup_t |a^R(t)-b^R(t)| = \sup_t |a^L(t)-b^L(t)| \le \sup_t |a(t)-b(t)| \,.
	\]
	Indeed, for every $t$,
	\[
		|a^R(t)-b^R(t)| = \lim_{s\downarrow t}|a(s)-b(s)|
		\le \sup_{r\in\mathbb R}|a(r)-b(r)| \,.
	\]
	Taking the supremum over $t$ gives
	\[
		\sup_t |a^R(t)-b^R(t)| \le \sup_t |a(t)-b(t)| \,.
	\]
	The same argument from the left gives
	\[
		\sup_t |a^L(t)-b^L(t)| \le \sup_t |a(t)-b(t)| \,.
	\]
	To compare the two versions, note that for every $h\in \mathsf{M}$
	\[
		h^R=[h^L]^R,\qquad  h^L=[h^R]^L \,.
	\]
	Applying the right-continuous inequality above to $a^L$ and $b^L$ gives
	\[
		\sup_t |a^R(t)-b^R(t)| = \sup_t |(a^L)^R(t)-(b^L)^R(t)|
		\le \sup_t |a^L(t)-b^L(t)| \,,
	\]
	and
	\[
		\sup_t |a^L(t)-b^L(t)| = \sup_t |(a^R)^L(t)-(b^R)^L(t)|
		\le \sup_t |a^R(t)-b^R(t)| \,.
	\]
	Thus,
	\[
		\sup_t |a^R(t)-b^R(t)| = \sup_t |a^L(t)-b^L(t)| \,.
	\]

	Fix $i$, $z\in\mathcal Z$, and $z_i^\ast\in\mathcal Z_i$. Applying the
	deterministic fact with
	\[
		a(t):=f(z_i^\ast,z_{-i},t), \qquad b(t):=f(z,t)
	\]
	gives
	\begin{align*}
		\sup_t \left| f^R(z_i^\ast,z_{-i},t)-f^R(z,t) \right|
		 & = \sup_t \left| f^L(z_i^\ast,z_{-i},t)-f^L(z,t) \right|   \\
		 & \le \sup_t \left| f(z_i^\ast,z_{-i},t)-f(z,t) \right| \,.
	\end{align*}
	Taking the supremum over $z$ and $z_i^\ast$ gives
	\[
		c_i(f^R)=c_i(f^L)\le c_i(f) \,.
	\]

	Dominated convergence, applied at each fixed $t$, gives
	\[
		\Exf{f^R_Z(t)} = \lim_{s\downarrow t}\Exf{f_Z(s)}, \qquad
		\Exf{f^L_Z(t)} = \lim_{s\uparrow t}\Exf{f_Z(s)} \,.
	\]
	Denoting by $\Exf{g_Z}$ the function $t\mapsto \Exf{g_Z(t)}$, we obtain
	\[
		\Exf{f^R_Z} = (\Exf{f_Z})^R, \qquad \Exf{f^L_Z} = (\Exf{f_Z})^L \,.
	\]
	The functions $\Exf{f^R_Z}$, $\Exf{f^L_Z}$, and $\Exf{f_Z}$ belong to
	$\mathsf{M}$.

	So fix $z\in\mathcal Z$.  Applying
	\[
		a(t):=f_z(t), \qquad b(t):=\Exf{f_Z(t)} \,,
	\]
	gives
	\[
		\sup_{t\in\mathbb R} \left| f^R_z(t)-\Exf{f^R_Z(t)} \right|
		= \sup_{t\in\mathbb R} \left| f^L_z(t)-\Exf{f^L_Z(t)} \right|
		\le \sup_{t\in\mathbb R} \left| f_z(t)-\Exf{f_Z(t)} \right| \,.
	\]
	This completes the proof.
\end{proof}

\subsection{Monotonicity}

\begin{lemma}\label{lem:interpolation_monotonicity}
	Let $g_Z$ be the interpolated function constructed in the proof
	of Lemma~\ref{lem:uniformisation}. For each fixed $z \in \mathcal{Z}$, the
	mapping $\vartheta \mapsto g(z, \vartheta)$ is non-decreasing.
\end{lemma}

\begin{proof}
	Let $\vartheta_1 < \vartheta_2$. Because the expectation $\Exf{f^{(2)}_Z(t)}$
	is strictly increasing, the generalised inverse gives
	$t_{\vartheta_1} \le t_{\vartheta_2}$. Two cases arise:
	\begin{itemize}
		\item If $t_{\vartheta_1} < t_{\vartheta_2}$, the monotonicity of
		      $t \mapsto f^{(2)}(z,t)$ bounds the maximum possible interpolated value at
		      $t_{\vartheta_1}$ by the minimum possible interpolated value at
		      $t_{\vartheta_2}$. Thus,
		      \[
			      g(z, \vartheta_1) \le f^{(2)\mathrm{R}}(z, t_{\vartheta_1})
			      \le f^{(2)\mathrm{L}}(z, t_{\vartheta_2}) \le g(z, \vartheta_2) \,.
		      \]
		\item If $t_{\vartheta_1} = t_{\vartheta_2} = t^\ast$, the interpolation
		      occurs across a single vertical jump. Because $\vartheta_1 < \vartheta_2$,
		      the unique interpolation weights for $\vartheta_2$ shift mass from the
		      left limit $f^{(2)\mathrm{L}}(z, t^\ast)$ towards the evaluation
		      $f^{(2)}(z, t^\ast)$ and the right limit $f^{(2)\mathrm{R}}(z, t^\ast)$.
		      The pointwise ordering
		      $f^{(2)\mathrm{L}}(z, t^\ast) \le f^{(2)}(z, t^\ast) \le
			      f^{(2)\mathrm{R}}(z, t^\ast)$ ensures this rightward weight shift gives
		      $g(z, \vartheta_1) \le g(z, \vartheta_2)$.
	\end{itemize}
	Thus, $g_Z$ is non-decreasing.
\end{proof}

\section{The upper bound}

\subsection{Lattice properties of the effective sample size}

Let $f \colon \mathcal{Z} \times \Theta \to \mathbb{R}$ be bounded, and
define the \emph{lattice functions}
$\overline{f^+}, \overline{f^-}, \overline{|f|} \colon \mathcal{Z} \to \mathbb{R}$ by
\[
	\overline{f^+}(z) = \sup_{\theta\in\Theta} (f(z,\theta))^+, \qquad
	\overline{f^-}(z) = \sup_{\theta\in\Theta} (f(z,\theta))^-, \qquad
	\overline{|f|}(z) = \sup_{\theta\in\Theta} |f(z,\theta)| \,.
\]
Each is a bounded function on $\mathcal{Z}$, hence a function on
$\mathcal{Z}\times\{\ast\}$ with singleton index set, and so has an
effective sample size.

\begin{lemma}\label{lem:lattice}
	The lattice functions have effective sample size at least that of $f$; that is,
	\[
		\eta_{\overline{f^+}} \ge \eta_f, \qquad
		\eta_{\overline{f^-}} \ge \eta_f, \qquad
		\eta_{\overline{|f|}} \ge \eta_f \,.
	\]
\end{lemma}

\begin{proof}
	For each coordinate $i$,
	\[
		c_i(\overline{f^+}) \le c_i(f), \qquad
		c_i(\overline{f^-}) \le c_i(f), \qquad
		c_i(\overline{|f|}) \le c_i(f)\,,
	\]
	and the result follows.
\end{proof}

We use the following in the proof of Lemma~\ref{prop:upper-envelope} for
the upper bound of the main theorem.

\begin{corollary}\label{cor:lattice_mcdiarmid}
	Let $f_Z \colon \Theta \to \mathbb{R}$ be a randomised function with finite
	$\eta_f$, and let $\ds \colon \Theta \to \mathbb{R}$ be a bounded deterministic
	function. If
	\[
		X(Z) := \sup_{\theta\in \Theta} |f_Z(\theta) - \ds(\theta)|
	\]
	is a measurable random variable, then
	\[
		\mathbb{P}\left\{ \sqrt{\eta_f}\left( X(Z) - \Exf{X(Z)} \right) \ge
		\varepsilon \right\} \le \e^{-2\varepsilon^2}, \qquad
		\text{for all } \varepsilon > 0 \,.
	\]
\end{corollary}

\begin{proof}
	Write $g(z,\theta) = f(z,\theta) - \ds(\theta)$. Because $\ds$ is deterministic,
	$c_i(g) = c_i(f)$ for all $i$, so $\eta_g = \eta_f$.
	Lemma~\ref{lem:lattice} gives $\eta_X \ge \eta_g = \eta_f$.
	If $\eta_X = \infty$, then $X(Z) = \Exf{X(Z)}$ surely and the bound holds.
	Otherwise $\infty > \eta_X \ge \eta_f > 0$, and for all $\varepsilon > 0$
	McDiarmid's inequality gives
	\[
		\mathbb{P}\left\{ X(Z) - \Exf{X(Z)} \ge \varepsilon \right\}
		\le \e^{-2\varepsilon^2\eta_X}
		\le \e^{-2\varepsilon^2\eta_f}\,.
	\]
	With $\varepsilon/\sqrt{\eta_f}$ in place of $\varepsilon$, the bound follows.
\end{proof}

\begin{remark}\label{rem:lattice}
	Consider the pointwise coefficient
	\[
		\frac{1}{\gamma_f} := \sup_{\theta\in \Theta}\sum_{i=1}^n
		c_i(f, \theta)^2, \qquad\text{where}\qquad
		c_i(f, \theta) := \sup_{(z^\ast_i,z)\in \mathcal{Z}_i\times \mathcal{Z}}
		\left| f_\theta(z_i^\ast,z_{-i})-f_\theta(z) \right| \,.
	\]
	Because the supremum over $\Theta$ is taken after summing over
	coordinates, rather than separately in each coordinate, we have
	$\gamma_f\ge\eta_f$. The pointwise McDiarmid--Hoeffding bounds
	\eqref{eq:macd}, therefore, hold with $\gamma_f$ in place of $\eta_f$,
	giving a pointwise sharper coefficient. However, $\gamma_f$ does not
	satisfy the lattice inequalities of Lemma~\ref{lem:lattice}, which are
	needed in the proof of Corollary~\ref{cor:lattice_mcdiarmid}.
\end{remark}

\subsection{The ramp property}

\begin{proof}[Proof of Lemma~\ref{lem:ramp-sup}]
	\emph{Positive deviation.} Let
	$S := \sup_{t\in[0,1]} (f(t)-t)^+ \ge \varepsilon$.
	Choose a sequence $(t_n) \subset [0,1]$ such that $f(t_n)-t_n \to S$.
	By compactness, passing to a subsequence if necessary, assume
	$t_n \to t^\ast \in [0,1]$.
	Because $f \le 1$, $t_n \le 1 - (f(t_n) - t_n)$. Taking limits gives
	$t^\ast \le 1 - S \le 1 - \varepsilon$.
	Thus, $t^\ast + \alpha \in [0,1]$ for all $\alpha \in [0, \varepsilon]$.

	Fix $\alpha \in (0, \varepsilon]$. For sufficiently large $n$,
	$t_n \le t^\ast + \alpha$. Monotonicity of $f$ implies
	$f(t^\ast + \alpha) \ge f(t_n)$. Thus,
	\begin{align*}
		f(t^\ast + \alpha) - (t^\ast + \alpha)
		 & \ge f(t_n) - t^\ast - \alpha                 \\
		 & = (f(t_n) - t_n) + t_n - t^\ast - \alpha \,.
	\end{align*}
	Letting $n \to \infty$ gives
	$f(t^\ast + \alpha) - (t^\ast + \alpha) \ge S - \alpha \ge
		\varepsilon - \alpha \ge 0$.

	\medskip

	\emph{Negative deviation.} Let
	$x := \sup_{t\in[0,1]} (t-f(t))^+ \ge \varepsilon$.
	Choose $(t_n) \subset [0,1]$ such that $t_n - f(t_n) \to x$, and assume
	$t_n \to t^\ast \in [0,1]$.
	Because $f \ge 0$, $t_n \ge t_n - f(t_n)$, giving
	$t^\ast \ge x \ge \varepsilon$.
	Thus, $t^\ast - \alpha \in [0,1]$ for all $\alpha \in (0, \varepsilon]$.

	Fix $\alpha \in (0, \varepsilon]$. For large $n$, $t_n \ge t^\ast - \alpha$.
	Monotonicity of $f$ gives $f(t^\ast - \alpha) \le f(t_n)$, and therefore
	\begin{align*}
		(t^\ast - \alpha) - f(t^\ast - \alpha)
		 & \ge t^\ast - \alpha - f(t_n)                 \\
		 & = (t_n - f(t_n)) + t^\ast - t_n - \alpha \,.
	\end{align*}
	Letting $n \to \infty$ gives us
	$(t^\ast - \alpha) - f(t^\ast - \alpha) \ge x - \alpha \ge
		\varepsilon - \alpha \ge 0$.
\end{proof}

\begin{proof}[Proof of Lemma~\ref{lem:scalar-envelope}]

	Write $c:=2\ln2$, and set
	\[
		\Phi(y):=\ln\left( 4\sqrt{2\pi}\,y+4\sqrt{y}\,\e^{-y/2} +\e^{-y} \right), \qquad \w(y)=\frac{\Phi(y)}{y} \,.
	\]
	Recall that $y_0$ solves
	\[
		\ln(\e+y)=\sqrt{y} \,.
	\]
	Numerically, $y_0>3.1028$.

	\medskip
	\emph{Reduction by squaring.}
	Multiplying the claimed bound by $\sqrt{c/y}>0$ and isolating the
	square root shows that it is equivalent to
	\[
		\sqrt{\left(1+\frac{\Phi(y)}{y}\right)^2+\frac cy}
		\le 1+\frac{2\sqrt{c}\,\ln(\e+y)}{y} \,.
	\]
	Both sides are positive. Squaring and multiplying by $y^2$ gives the
	equivalent inequality
	\[
		2y\Phi(y)+\Phi(y)^2+cy \le 4\sqrt{c}\,y\ln(\e+y) +4c\,\ln(\e+y)^2 \,.
	\]

	\medskip
	\emph{Majorising $\Phi$.}
	The function
	\[
		y\mapsto 4\sqrt{y}\,\e^{-y/2}+\e^{-y}
	\]
	is decreasing for $y\ge1$. Hence, for $y\ge3$,
	\[
		4\sqrt{y}\,\e^{-y/2}+\e^{-y} \le 4\sqrt{3}\,\e^{-3/2}+\e^{-3} \le (11-4\sqrt{2\pi})\,3 \,.
	\]
	Because $y\ge y_0>3$, we have
	\[
		4\sqrt{2\pi}\, y +4\sqrt{y}\,\e^{-y/2} +\e^{-y}
		\le 11y \,.
	\]
	Therefore
	\[
		0<\Phi(y)\le \ln(11y), \qquad y\ge y_0 \,.
	\]
	Because the map $\xi\mapsto2y\xi+\xi^2$ is increasing on $\xi\ge0$,
	it is enough to prove
	\[
		2y\ln(11y)+\ln(11y)^2+cy \le 4\sqrt c\,y\ln(\e+y) +4c\,\ln(\e+y)^2 \,.
	\]
	Define
	\[
		\Psi(y) := 4\sqrt c\,y\ln(\e+y) +4c\,\ln(\e+y)^2 -2y\ln(11y)-\ln(11y)^2-cy \,.
	\]
	We prove that $\Psi(y)\ge0$ for $y\ge y_0$.

	\medskip
	\emph{Convexity of $\Psi$.}
	Differentiating,
	\[
		\Psi'(y) = 4\sqrt c\left(\ln(\e+y)+\frac{y}{\e+y}\right) +\frac{8c\ln(\e+y)}{\e+y} -2\ln(11y)-2-\frac{2\ln(11y)}{y}-c
	\]
	and
	\[
		\Psi''(y) = \frac{4\sqrt{c}}{\e+y} +\frac{4\sqrt c\,\e}{(\e+y)^2}
		+\frac{8c(1-\ln(\e+y))}{(\e+y)^2}
		-\frac2y +\frac{2(\ln(11y)-1)}{y^2} \,.
	\]
	We show that $\Psi''(y)\ge0$ for $y\ge y_0$.

	The only negative terms in $\Psi''$ are $-2/y$ and the negative part of
	the term containing $1-\ln(\e+y)$. We use
	\[
		\frac{8c(\ln(\e+y)-1)}{(\e+y)^2}
		\le \frac{4c}{\e^3} \,,
	\]
	because
	\[
		\max_{u>0}\frac{\ln(u)-1}{u^2}=\frac{1}{2\e^3} \,.
	\]

	\medskip
	\emph{The tail ($y\ge 6$).}
	For $y\ge 6$, discard the non-negative terms
	\[
		\frac{4\sqrt c\,\e}{(\e+y)^2}, \qquad
		\frac{2(\ln(11y)-1)}{y^2}
	\]
	from $\Psi''$. After multiplying by $(\e+y)^2$, it is enough to
	prove
	\[
		\left(\frac{4\sqrt c}{\e+y}-\frac2y\right)
		(\e+y)^2 -8c(\ln(\e+y)-1) \ge 0 \,.
	\]
	The first term equals
	\[
		(4\sqrt c-2)y +(4\sqrt c-4)\e -\frac{2\e^2}{y} \,.
	\]
	For $y\ge6$,
	\[
		\frac{2\e^2}{y}\le\frac{\e^2}{3} \,.
	\]
	Thus, it is enough to prove $\Xi(y)\ge0$, where
	\[
		\Xi(y) := (4\sqrt c-2)y +(4\sqrt c-4)
		\e -\frac{\e^2}{3} - 8c(\ln(\e+y)-1) \,.
	\]
	But
	\[
		\Xi'(y)
		= (4\sqrt c-2)-\frac{8c}{\e+y}>0, \qquad y\ge6
	\]
	and
	\[
		\Xi(6)>2.79 \,.
	\]
	Hence, $\Psi''(y)\ge0$ for $y\ge6$.

	\medskip
	\emph{The bounded range ($y_0\le y\le6$).}
	Using the global bound above in the formula for $\Psi''$, we have
	\[
		\Psi''(y) \ge \frac{4\sqrt c}{\e+y} +\frac{4\sqrt c\,\e}{(\e+y)^2} +\frac{2(\ln(11y)-1)}{y^2} -\frac2y-\frac{4c}{\e^3} \,.
	\]
	The function
	\[
		y\mapsto\frac{\ln(11y)-1}{y^2}
	\]
	is decreasing on $[y_0,\infty)$, because
	\[
		\frac{\mathrm d}{\mathrm dy} \left(\frac{\ln(11y)-1}{y^2}\right) = \frac{3-2\ln(11y)}{y^3}<0 \,.
	\]
	Hence, on $[y_0,4]$,
	\[
		\Psi''(y) \ge \frac{4\sqrt c}{4+\e} +\frac{4\sqrt c\,\e}{(4+\e)^2} +\frac{2(\ln44-1)}{16} -\frac{2}{y_0} -\frac{4c}{\e^3} >0 \,.
	\]
	Similarly, on $[4,6]$,
	\[
		\Psi''(y) \ge \frac{4\sqrt c}{6+\e} +\frac{4\sqrt c\,\e}{(6+\e)^2} +\frac{2(\ln44-1)}{36} -\frac12 -\frac{4c}{\e^3} >0 \,.
	\]
	Therefore, $\Psi''(y)\ge0$ for all $y\ge y_0$.

	\medskip 
	Because direct verification shows
	\[
		\Psi'(y_0)>0 \,,
	\]
	and
	\[
		\Psi(y_0)>0 \,,
	\]
	the function $\Psi$ is increasing and positive on $[y_0,\infty)$. This
	proves the majorised inequality; hence, the squared inequality, and
	therefore, the scalar envelope bound.
\end{proof}

\bibliographystyle{imsart-number}
\bibliography{prob}

@article{kelley1959,
  author = {Kelley, J. L.},
  title = {Measures on {B}oolean algebras},
  journal = {Pacific Journal of Mathematics},
  volume = {9},
  number = {4},
  pages = {1165--1177},
  year = {1959},
}

@book{boucheron2013concentration,
  author = {Boucheron, St{\'e}phane and Lugosi, G{\'a}bor and Massart, Pascal},
  title = {Concentration inequalities: {A} nonasymptotic theory of independence},
  publisher = {Oxford University Press},
  year = {2013},
}

@article{canay2018wild,
  author = {Canay, Ivan A. and Santos, Andres and Shaikh, Azeem M.},
  title = {The wild bootstrap with a ``small'' number of ``large'' clusters},
  journal = {The Review of Economics and Statistics},
  volume = {103},
  number = {2},
  pages = {346--363},
  year = {2021},
  doi = {10.1162/rest_a_00887},
}

@article{debreu1964continuity,
  author = {Debreu, Gerard},
  title = {Continuity properties of {P}aretian utility},
  journal = {International Economic Review},
  volume = {5},
  pages = {285--293},
  year = {1964},
}

@article{dibenedetto2021nonexchangeable,
  author = {Di Benedetto, Giuseppe and Caron, Fran{\c{c}}ois and Teh, Yee Whye},
  title = {Nonexchangeable random partition models for microclustering},
  journal = {The Annals of Statistics},
  volume = {49},
  pages = {1931--1957},
  year = {2021},
  doi = {10.1214/20-AOS2003},
}

@article{diebolt1996density,
  author = {Diebolt, Jean and Posse, Christian},
  title = {On the density of the maximum of smooth {G}aussian processes},
  journal = {The Annals of Probability},
  volume = {24},
  pages = {1104--1129},
  year = {1996},
}

@article{hoeffding1963probability,
  author = {Hoeffding, Wassily},
  title = {Probability inequalities for sums of bounded random variables},
  journal = {Journal of the American Statistical Association},
  volume = {58},
  pages = {13--30},
  year = {1963},
}

@article{kingman1978representation,
  author = {Kingman, John F. C.},
  title = {The representation of partition structures},
  journal = {Journal of the London Mathematical Society},
  volume = {18},
  pages = {374--380},
  year = {1978},
}

@book{ledoux2001concentration,
  author = {Ledoux, Michel},
  title = {The concentration of measure phenomenon},
  series = {Mathematical Surveys and Monographs},
  volume = {89},
  publisher = {American Mathematical Society},
  year = {2001},
}

@book{macwilliams1977theory,
  author = {MacWilliams, Florence Jessie and Sloane, Neil James Alexander},
  title = {The theory of error-correcting codes},
  series = {North-Holland Mathematical Library},
  volume = {16},
  publisher = {North-Holland},
  address = {Amsterdam},
  year = {1977},
}

@article{massart1990tight,
  author = {Massart, Pascal},
  title = {The tight constant in the {D}voretzky--{K}iefer--{W}olfowitz
           inequality },
  journal = {The Annals of Probability},
  volume = {18},
  pages = {1269--1283},
  year = {1990},
}

@incollection{mcdiarmid_1989,
  author = {McDiarmid, Colin},
  title = {On the method of bounded differences},
  booktitle = {Surveys in Combinatorics, 1989},
  series = {London Mathematical Society Lecture Note Series},
  volume = {141},
  pages = {148--188},
  publisher = {Cambridge University Press},
  year = {1989},
}

@article{pitman1995exchangeable,
  author = {Pitman, Jim},
  title = {Exchangeable and partially exchangeable random partitions},
  journal = {Probability Theory and Related Fields},
  volume = {102},
  pages = {145--158},
  year = {1995},
}

@article{talagrand1996new,
  author = {Talagrand, Michel},
  title = {A new look at independence},
  journal = {The Annals of Probability},
  volume = {24},
  pages = {1--34},
  year = {1996},
}

@book{vanderVaartWellner1996,
  author = {van der Vaart, Aad W. and Wellner, Jon A.},
  title = {Weak convergence and empirical processes: {W}ith applications to
           statistics},
  series = {Springer Series in Statistics},
  publisher = {Springer},
  address = {New York},
  year = {1996},
  doi = {10.1007/978-1-4757-2545-2},
}

@book{vershynin2018high,
  author = {Vershynin, Roman},
  title = {High-dimensional probability: {A}n introduction with applications in
           data science},
  series = {Cambridge Series in Statistical and Probabilistic Mathematics},
  publisher = {Cambridge University Press},
  year = {2018},
}

@book{wainwright2019high,
  author = {Wainwright, Martin J.},
  title = {High-dimensional statistics: {A} non-asymptotic viewpoint},
  series = {Cambridge Series in Statistical and Probabilistic Mathematics},
  volume = {48},
  publisher = {Cambridge University Press},
  year = {2019},
}

\end{document}